\documentclass[12pt,a4paper,leqno,verbatim]{amsart}
\newcounter{minutes}
\divide\time by 60
\newcounter{hours}
\multiply\time by 60 \addtocounter{minutes}{-\time}

\usepackage{amssymb}
\usepackage{hyperref}
\usepackage{graphicx}
\usepackage{color}

\date{}
\newfont{\cyrilic}{wncyr10 scaled 1000}
\title[On Jordan's Redheffer's and Wilker's inequality]{On Jordan's Redheffer's and Wilker's inequality}

\author[B.A Bhayo]{Barkat Ali Bhayo}

\author[J. S\'andor]{J\'ozsef S\'andor}

\newcommand{\comment}[1]{}

\swapnumbers
\theoremstyle{plain}

\newtheorem{theorem}[equation]{Theorem}
\newtheorem{lemma}[equation]{Lemma}

\newtheorem{corollary}[equation]{Corollary}

\newtheorem{remark}[equation]{Remark}

\numberwithin{equation}{section}

\begin{document}




\begin{abstract}
In this paper, the authors offer new Jordan, Redheffer and Wilker type
inequalities, along with refinements and converses. Connections with
Euler's gamma function are pointed out, too.

\end{abstract}

\maketitle

%
%
%
%


\section{Introduction}

In the recent years, the refinements of the inequalities of trigonometric functions such as
Adamovi\'c-Mitrinovi\'c, Cusa-Huygens, Jordan inequality, Redheffer inequality,
Becker-Stark inequality, Wilker inequality,
Huygens inequality, and Kober inequality
have been studied extremely by numerous authors,
e.g., see \cite{avv1,avv2,czq,kvv,lili,neua,neu2a,neusana,nejoz1a,joz111a,joz1a,joz2,joz3,sanben}
and the references therein. Motivated by these rapid studies,
in this paper we make a contribution
to the subject by refining the Cusa-Huygens, Jordan and Redheffer inequality, and establish a Wilker type inequality.
In all cases, we
give the upper and lower
bound of $\sin(x)/x$ in terms of simple functions. Meanwhile, we give some Redheffer type inequalities for
trigonometric functions, which refine the existing results in the literature.

For the representation of trigonometric functions in terms of gamma function, we define
for $x,y>0$  the classical \emph{gamma function} $\Gamma$, the {\it digamma function} $\psi$ and the \emph{beta function}  $B(\cdot,\cdot)$ by
   $$\Gamma(x) = \int^\infty_0 e^{-t}t^{x-1}\,dt,\quad
       \psi(x) = \frac{\Gamma'(x)}{\Gamma(x)},\quad
        B(x,y) = \frac{\Gamma(x)\Gamma(y)}{\Gamma(x+y)},$$
respectively. We denote the $trigamma$ function $\psi'$ by $K$, and defined as,
\begin{equation}\label{trigama}
K(x)=\sum_{n=0}^\infty\frac{1}{(x+n)^2},\quad K'(x)=-2\sum_{n=0}^\infty\frac{1}{(x+n)^3}.
\end{equation}
We refer to reader to see \cite{joz0,joz1,joz2aa,joz3aa} for more properties and relations of $K$.
The functions $\Gamma$ and $\psi$ satisfy the recurrence relation
\begin{equation}\label{recgam}
\Gamma(1+z)=z\Gamma(z),\quad \psi(1+z)=z\psi(z).
\end{equation}
The following relation
\begin{equation}\label{recpol}
\psi(1+z)-\psi(z)=\frac{1}{z}
\end{equation}
follows from \eqref{recgam}. Differentiating both sides of \eqref{recpol}, one has
\begin{equation}\label{recK}
K(1+z)-K(z)=-\frac{1}{z^2}.
\end{equation}

The following Euler's reflection formula \cite[6.1.17]{as}
\begin{equation}\label{eulref}
\Gamma(t)\Gamma(1-t)=\frac{\pi}{\sin (\pi t)},\quad 0<t<1,
\end{equation}
can be written as
\begin{equation}\label{sinref}
\frac{x}{\sin(x)}=\Gamma\left(1+\frac{x}{\pi}\right)\Gamma\left(1-\frac{x}{\pi}\right)
=B\left(1+\frac{x}{\pi},1-\frac{x}{\pi}\right),\quad 0<x<\pi.
\end{equation}
The logarithmic differentiation to both sides of \eqref{eulref} gives the following reflection formula,

\begin{equation}\label{eulreftan}
\psi(1-t)-\psi(t)=\frac{\pi}{\tan (\pi t)}.
\end{equation}

Replacing $t$ by $t+1/2$ in \eqref{eulref}, we get
\begin{equation}\label{cosref}
\frac{x}{\cos(x)}=\frac{x}{\pi}\Gamma\left(\frac{1}{2}+\frac{x}{\pi}\right)\Gamma\left(\frac{1}{2}-\frac{x}{\pi}\right)
,\quad 0<x<\frac{\pi}{2}.
\end{equation}


Next we recall the Adamovi\'c-Mitrinovi\'c \cite[p.238]{mit} and Cusa-Huygens \cite{sanben} inqualities
\begin{equation}\label{laztri}
(\cos x)^{1/3}<\frac{\sin x}{x}<\frac{\cos x+2}{3},\quad 0<|x|<\frac{\pi}{2},
\end{equation}
For the refinement of \eqref{laztri}, see \cite{kvv,neua,neusana,joz3,sanben,yang}.

Our first main result refines the inequalities in \eqref{laztri} as follows:

\begin{theorem}\label{thm0803} For $x\in(0,\pi)$, the following inequalities hold,
\begin{equation}\label{inequ0703}
\frac{2-(x/\pi)^3+\cos(x)}{3}<\frac{\sin(x)}{x}<\frac{2-(x/\pi)^4+\cos(x)}{3}.
\end{equation}
\end{theorem}

%

In \cite{wilker}, Wilker asked a question to find the largest constant $c$ such that the inequality
\begin{equation}\label{wilker1003}
a(x)>c\,b(x),\quad 0<|x|<\frac{\pi}{2},
\end{equation}
holds true,
where $$a(x)=\left(\frac{\sin(x)}{x}\right)^2+\frac{\tan(x)}{x}-2,\quad {\rm and}\quad b(x)=x^3\tan(x).$$
The inequality \eqref{wilker1003} is known as the Wilker inequality in the literature. 
Anglesio \cite{sumner}
proved that the ratio $a(x)/b(x)$ is decreasing in $x\in(0,\pi/2)$,
and he answered the question by showing the following two sided
inequality,
\begin{equation}\label{angineq}
2+\frac{8}{45}x^3\tan(x)<\left(\frac{\sin(x)}{x}\right)^2+\frac{\tan(x)}{x}<2+\frac{16}{\pi^4}x^3\tan(x),
\end{equation}
with the best possible constants $8/45$ and $16/\pi^4$. For the new proofs  and refinement of \eqref{angineq},
 see \cite{guo,Pin,wang,wu2,zz,z3}. Thereafter the following Wilker type inequalities,
\begin{eqnarray*}
3+\frac{3}{20}x^3\tan(x)&<&2\,\frac{\sin(x)}{x}+\frac{\tan(x)}{x}
<3+\frac{16}{\pi^4}x^3\tan(x),\,0<x<\frac{\pi}{2},\\
2+\frac{2}{45}x^3\sin(x)&<&\left(\frac{x}{\sin(x)}\right)^2+\frac{x}{\tan(x)}
<2+\left(\frac{2}{\pi}-\frac{16}{\pi^3}\right)x^3\sin(x),\,0<x<\frac{\pi}{2},\\
& &\left(\frac{\sinh(x)}{x}\right)^2+\frac{\tanh(x)}{x}
>2+\frac{8}{45}x^3\tanh(x),\,x>0,\\
& &\left(\frac{x}{\sinh(x)}\right)^2+\frac{x}{\tanh(x)}
<2+\frac{2}{45}x^3\sinh(x),\,x>0,
\end{eqnarray*}
were established by Chen-S\'andor \cite{chen-s}, Wang \cite{wang}, Zhu \cite{z4}, and Sun-Zhu \cite{sun-z},
 respectively.

We establish an other Wilker type inequality by giving the following result.
\begin{theorem}\label{wu-sri} For $x\in(0,\pi)$, we have
$$3+\left(\frac{x}{\pi}\right)^4\frac{x}{\sin(x)}<2\frac{x}{\sin(x)}+\frac{x}{\tan(x)}
<3+\left(\frac{x}{\pi}\right)^3\frac{x}{\sin(x)}.$$
\end{theorem}

The following inequality
\begin{equation}\label{wusri}
2\frac{x}{\sin(x)}+\frac{x}{\tan(x)}>2,\quad 0<|x|<\frac{\pi}{2},
\end{equation}
has been recently established by Wu and Srivastava \cite{wu}, which is sometime known as the second Wilker inequality \cite{neusana}. Inequality \eqref{wusri} can be written as,
$$\frac{\sin(x)}{x}<\frac{1}{2}(\frac{x}{\sin(x)}+\cos(x)),\quad 0<|x|<\frac{\pi}{2}.$$
It was shown in \cite{neusana} that this inequality is weaker than the Cusa-Huygens inequality 
\eqref{laztri}, as follows:
$$\frac{\sin(x)}{x}<\frac{2+\cos(x)}{3}<\frac{1}{2}\left(\frac{x}{\sin(x)}+\cos(x)\right),\quad 0<|x|<\frac{\pi}{2},$$
here second inequality is equivalent to
\begin{equation}\label{ineq1103}
3\,\frac{x}{\sin(x)}+\cos(x)>4,\quad 0<|x|<\frac{\pi}{2}.
\end{equation}
Recently, Mortici \cite{mor} refined this inequality \eqref{ineq1103} as follows:
\begin{equation}\label{neuineq}
3\,\frac{x}{\sin(x)}+\cos(x)>4+\frac{1}{10}x^4+\frac{1}{210}x^6,\quad 0<|x|<\frac{\pi}{2}.
\end{equation}
We establish the following sharp result as a counterpart of \eqref{neuineq}.
\begin{theorem}\label{1003} For $x\in(0,\pi)$, we have
$$2-\frac{x^4}{\pi^4}<3\,\frac{\sin(x)}{x}-\cos(x)<2-\frac{x^3}{\pi^3}.$$
\end{theorem}

The well-know Jordan's inequality \cite{mit} states,
$$\frac{\pi}{2}\leq \frac{\sin(x)}{x},\quad 0<x\leq \frac{\pi}{2},$$
equality with $x=\pi/2$.

In \cite{deb}, Debnath and Zhao refined the Jordan's inequality as below,
\begin{equation}\label{debinequ1}
dz_l(x)=\frac{2}{\pi }+\frac{1}{12\pi}(\pi ^2-4 x^2)\leq\frac{\sin(x)}{x},
  \end{equation}
    \begin{equation}\label{debinequ2}
DZ_l(x)=\frac{2}{\pi }+\frac{1}{\pi^3}(\pi ^2-4 x^2)\leq\frac{\sin(x)}{x},
  \end{equation}
    for $x\in(0,\pi/2)$, equality in both inequalities with $x=\pi/2$.

    The following inequality
    \begin{equation}\label{ozbaninequ1}
 O_l(x)=\frac{2}{\pi }+\frac{1}{\pi^3}(\pi ^2-4 x^2)+\frac{4(\pi-3)}{\pi^3}\left(x-\frac{\pi}{2}\right)^2\leq
\frac{\sin(x)}{x},
  \end{equation}
     for $x\in(0,\pi/2)$ is due to \"Ozban \cite{ozban}, equality with $x=\pi/2$.

    In \cite{z2}, Zhu proved that
    \begin{equation}\label{zinequj}
\frac{\sin(x)}{x}\leq\frac{2}{\pi }+\frac{(\pi -2) \left(\pi ^2-4
   x^2\right)}{\pi^3}=Z_u(x),
\end{equation}
for $x\in(0,\pi/2)$.

For the following inequalities
\begin{equation}\label{doninequ}
J_l(x)=\frac{2}{\pi }+\frac{\pi ^4-16 x^4}{2\pi^5}<\frac{\sin(x)}{x}<\frac{2}{\pi }+\frac{(\pi -2) \left(\pi ^4-16
   x^4\right)}{\pi^5}=J_u(x),
\end{equation}
$x\in(0,\pi/2)$, see \cite{jiang}.

In \cite[Theorem 1.3]{kvv}, Kl\'en et al. proved that
\begin{equation}\label{kvvineq}
1-\frac{x^2}{6}<\frac{\sin(x)}{x}< 1-\frac{2x^2}{3\pi^2},
\end{equation}
for $x\in(-\pi/2,\pi/2)$. The following inequality refines the second inequality in \eqref{kvvineq},
\begin{equation}\label{liineq}
\frac{\sin(x)}{x}\leq \frac{1-x^2/\pi^2}{\sqrt{1+3(x/\pi)^4}},
\end{equation}
for $x\in(0,\pi)$, see \cite{lili}.

Our first main result reads as follows:

\begin{theorem}\label{thm2} For $x\in(0,\pi/2)$, we have
$$C_l(x)<\frac{\sin (x)}{x}<C_u(x),$$
where
$$C_l(x)=\left(1-\frac{x^2}{\pi ^2}\right)^{\pi ^2/6}\quad {\rm and}\quad C_u(x)=\left(1-\frac{x^2}{\pi
   ^2}\right)^{3/2}.$$ 
    \end{theorem}

The right side of the following theorem can not be compared with the corresponding side of Theorem \ref{thm2},
and the second inequality is weaker than the above result.

\begin{theorem}\label{thm1} For $x\in(0,\pi)$, we have
$$D_l(x)<\frac{\sin(x)}{x}<D_u(x),$$
where
$$ D_l(x)=1-\left(\frac{x}{\pi}\right)^2\left(2-\frac{x}{\pi}\right) \quad {\rm and}\quad
D_u(x)=1-\frac{1}{2}\left(\frac{x}{\pi}\right)^2\left(3-\left(\frac{x}{\pi}\right)^2\right).$$
\end{theorem}
It is not difficult to see that Theorem \ref{thm1} refines the inequalities \eqref{kvvineq} and \eqref{liineq}.
Obviously, one can see that $dz_l(x)<DZ_l(x)<O_l(x)$, $J_l(x)<DZ_l(x)$, and
$Z_u(x)<J_u(x)$ for $x\in(0,\pi/2)$. Now it is natural to compare our result with \eqref{ozbaninequ1} and
\eqref{zinequj}. For this purpose
we give the following inequalities by using the Mathematica Software\textsuperscript{\tiny{\textregistered}} \cite{ru},
\begin{eqnarray*}
O_l(x)&<&C_l(x), \quad x\in(0,1.19540),\\
O_l(x)&<&D_l(x), \quad x\in(0,0.92409),\\
C_u(x)&<&Z_u(x), \quad x\in(0,1.09447),\\
D_u(x)&<&Z_u(x), \quad x\in(0,0.95784).
\end{eqnarray*}
We see that our result refines \eqref{ozbaninequ1} and
\eqref{zinequj} in the given interval of $x$.

The following result is the consequence of Theorem \ref{thm1}.
\begin{theorem}\label{gammathm} For $y\in(0,1)$, we have
\begin{enumerate}
\item $\displaystyle B(x,y)<\frac{1}{xy}\frac{x+y}{1+xy},\quad 0<x<1 $,
\item $\displaystyle\frac{1}{xy}\frac{x+y}{1+xy}<B(x,y),\quad 1<x<\infty$.\\
\end{enumerate}
\end{theorem}
It is worth to mention that the part (1) of Theorem \ref{gammathm} recently appeared in \cite{ivady} as a
one of the main results. The proof of the claim is based on
\cite[Lemma 2.5]{ivady}, and the proof of the lemma is invalid.

The following Redheffer inequality \cite{re}
\begin{equation}\label{sinRed}
\frac{\pi^2-x^2}{\pi^2+x^2}\leq \frac{\sin x}{x},\quad 0<x\leq \pi,
\end{equation}
was proved by Williams \cite{wi}.
Chen et al. \cite{czq} obtained the three Redheffer-type inequalities
for $\cos x,\,\cosh x$ and $(\sinh x)/x$. Sun and Zhu \cite{lili,sz} proved the Redheffer-type two-sided inequalities for trigonometric and hyperbolic functions. The inequalities appeared in \cite{lili} were refined in \cite{z2},
and read as follows:

\begin{equation}\label{zhu28a}
\left(\frac{\pi^2-x^2}{\sqrt{\pi^4+3x^4}}\right)^{\pi^2/6}<\frac{\sin(x)}{x}
<\left(\frac{\pi^2-x^2}{\sqrt{\pi^4+3x^4}}\right),\quad 0<x<\pi,
\end{equation}

\begin{equation}\label{zhu28b}
\left(\frac{\pi^2-4x^2}{\sqrt{\pi^4+48x^4}}\right)^{\pi^2/6}<\cos(x)
<\left(\frac{\pi^2-4x^2}{\sqrt{\pi^4+48x^4}}\right)^{3/4},\quad 0<x<\frac{\pi}{2}.
\end{equation}

\begin{equation}\label{zhu28c}
\left(\frac{\sqrt{\pi^4+48x^4}}{\pi^2-4x^2}\right)^{1/2}<\frac{\tan(x)}{x}
<\left(\frac{\sqrt{\pi^4+48x^4}}{\pi^2-4x^2}\right)^{\pi^2/6},\quad 0<x<\frac{\pi}{2}.
\end{equation}
One can see easily that the first inequality in Theorem \ref{thm1} refines the \eqref{sinRed}.
Our following theorem refines \eqref{zhu28a}, as well as Cusa-Huygens inequality,
and also gives a new upper bound for the right side
of \eqref{laztri}.

\begin{theorem}\label{thm3} For $x\in(0,\pi/2)$ one has
\begin{equation}\label{thm3ineq}
\frac{\sin(x)}{x}<\frac{\pi^2-x^2}{\pi^2+\alpha x^2}<\frac{2+\cos (x)}{3}<\frac{\pi^2-x^2}{\pi^2+\beta x^2},
\end{equation}
with best possible constants $\alpha=\pi^2/6-1\approx 0.644934$ and $\beta=1/2$.
\end{theorem}

\begin{theorem}\label{thm4} The function
$$g(t)=\frac{\log((1+t^2)/(1-t^2))}{\log(1/\cos(\pi t/2))}$$
is strictly decreasing from $(0,1)$ on $(\alpha,\beta)$. In particular, for $x\in(0,\pi)$
$$\cos(x/2)^\alpha<\frac{\pi^2-x^2}{\pi^2+x^2}<\cos(x/2)^\beta,$$
with the best possible constants $\alpha=16/\pi^2\approx 1.62114$ and $\beta=1$.
\end{theorem}

In the following theorem we refine the inequalities given in \eqref{zhu28b}.

\begin{theorem}\label{thm1702} The following function
$$h(t)=\frac{1-4t^2}{t^2\cos(\pi t)}-\frac{1}{t^2}$$
is strictly increasing from $(0,1/2)$ onto $((\pi^2-8)/2,16/\pi-4)$. In particular, for $x\in(0,\pi/2)$
$$\frac{\pi^2-4x^2}{\pi^2+\alpha x^2}<\cos(x)<\frac{\pi^2-4x^2}{\pi^2+\beta x^2},$$
with the best possible constants $\alpha=16/\pi-4\approx 1.09296$ and $\beta=(\pi^2-8)/2\approx 0.934802$.
\end{theorem}

The paper is organized into three sections as follows. Section 1,
contains the introduction and the statements of our main results. In Section 2,
we give some lemmas, which will be used
in our proofs.
Section 3 is consists of the proofs of the main results and some corollaries.
\section{Preliminaries}

The following result is sometime called the Monotone l'H\^opital rule, which is due to
Anderson et al. \cite{avv1}.
\begin{lemma}\label{lem0avv}
For $-\infty<a<b<\infty$,
let $f,g:[a,b]\to \mathbb{R}$
be continuous on $[a,b]$, and be differentiable on
$(a,b)$. Let $g^{'}(x)\neq 0$
on $(a,b)$. If $f^{'}(x)/g^{'}(x)$ is increasing
(decreasing) on $(a,b)$, then so are
$$\frac{f(x)-f(a)}{g(x)-g(a)}\quad and \quad \frac{f(x)-f(b)}{g(x)-g(b)}.$$
If $f^{'}(x)/g^{'}(x)$ is strictly monotone,
then the monotonicity in the conclusion
is also strict.
\end{lemma}

\begin{lemma}\label{lem0} For $a\in(0,1)$ and $x>0$, the following inequality holds,
$$\frac{1-a}{x+a}<\psi(1+x)-\psi(x+a).$$
\end{lemma}

\begin{proof} It is well-known that the function $f(x)= x\psi(x),\,x>0$ is strictly convex \cite[Theorem 6]{joz3aa}.
This implies that,
$$f(a\,r+(1-a)s)<a\,f(r)+(1-a)f(s),\quad r,s>0,\,a\in(0,1).$$
Setting $r=1+x$ and $s=x$ in the above inequality, we get
$$(a+x)\psi(a+x)<(a+ax)\psi(1+x)+(1-a)x\psi(x)),$$
now the proof follows easily if we replace $\psi(x)=\psi(1+x)-1/x$.
\end{proof}


\begin{lemma}\label{lem1} For $x\in(1,\infty)$ and $y\in(0,1)$, we have
$$\frac{\Gamma(1+x)\Gamma(1+y)}{\Gamma(1+x+y)}>\frac{1}{1+xy}.$$
\end{lemma}

\begin{proof} Let
$$f_y(x)=\log(\Gamma(1+x))+\log(\Gamma(1+y))-\log(\Gamma(1+x+y))+\log((1+xy)),$$
for $x\in(1,\infty)$ and $y\in(0,1)$, clearly $f_y(1)=0$.
Differentiating $f$ with respect to $x$ and using the formula \eqref{recpol}, we get
\begin{eqnarray*}
f'_y(x)&=&\frac{y}{1+xy}+\psi(1+x)-\psi(1+x+y)\\
&=&\frac{y}{1+xy}-\frac{1}{x+y}+\psi(1+x)-\psi(x+y)\\
&>&\frac{y}{1+xy}-\frac{1}{x+y}+\frac{1-y}{x+y}=\frac{y(1-y)(x-1)}{(x+y)(1+xy)}>0,
\end{eqnarray*}
by Lemma \ref{lem0}.
Thus, $f$ is strictly increasing, and $f_y(x)>f_y(1)=0$, this implies the proof of part (1).
For the proof of part (2), see \cite[(3.2)]{ivady}.
\end{proof}


\begin{remark} \rm
It is easy to see that the convexity of the function $x\mapsto \log(x\Gamma(x)),\,x>0,$
implies the following inequality,
\begin{equation}\label{gru1}
\frac{1}{xy}\left(\frac{x+y}{2}\right)^2<\frac{\Gamma(x)\Gamma(y)}{\Gamma((x+y)/2)^2}.
\end{equation}
Replacing $x$ by $1-x/\pi$ and $y$ by $1-x/\pi$ in \eqref{gru1}, and applying
\eqref{eulref}, we get
\begin{equation}\label{grcon}
\frac{\sin(x)}{x}<\frac{\pi^2-x^2}{\pi^2},\quad 0<x<\pi.
\end{equation}
This improves the  following inequality
$$\frac{\sin(x)}{x}\leq \left(\frac{\pi^2-x^2}{\pi^2+x^2}\right)^{1/2},\quad 0<x\leq \pi,$$
which was proved in \cite{joz01}.
\end{remark}

%

\section{Proof of the main results}

In this section we will give the proofs of the main results
highlighted in the first section, as well as some corollaries are
being established.

\vspace{.3cm}
\noindent{\bf Proof of Theorem \ref{thm0803}.} \rm The proof of the first inequality is trivial.
For the proof of the second and the third inequality, we define
$$f(x)=\frac{2+\cos(x)-3\sin(x)/x}{x^4},$$
for $x\in(0,\pi)$.
We will prove that $f$  is strictly decreasing from $(0, \pi)$ onto $(1/\pi^4, 1/60)$. One has,
$$x^6f(x)=-x^2\sin(x)-7x\cos(x)+15\sin(x)-8x=f_1(x),$$
$$f_1'(x)=5x\sin(x)+8\cos(x)-x^2\cos(x)-8,$$
$$f_2''(x)=-3\sin(x) +3x\cos(x) +x^2\sin (x),\, f_1'''(x)=x(x\cos(x)-\sin(x))<0.$$
Thus $f_1''(x)<f_1''(0)=0$, $f_1'(x)<f_1'(0)=0,$ and $f_1(x)<f_1(0)=0.$ The limiting values can be achieved by l'H\^opital rule.  $\hfill\square$

Similarly, for the proof of the first inequality we will prove that the function
$$g(x)=\frac{2+\cos(x)-3\sin(x)/x}{x^3}$$
is strictly increasing from $(0,\pi)$ onto $(0,1/\pi^3)$.
One has
\begin{eqnarray*}
x^5g'(x)&=& -(x^2)\sin(x) -6x\cos(x)+12 \sin(x)-6x=h(x),\\
h'(x)&= & 4x\sin(x)-x^2\cos(x)+6\cos(x)-6,\\
h''(x)&=& 2x\cos(x) +x^2\sin(x) -2\sin(x), \,h'''(x)= x^2\cos(x).
\end{eqnarray*}
Thus $h'''(x)$ is positive in $x\in(0, \pi/2)$ and negative in $x\in(\pi/2, \pi)$.  Since $h''(\pi/2)=\pi^2/4-2 >0$
 and $h''(\pi)=-2\pi <0$, and $h''(x)$ is strictly increasing (resp. decreasing) in $x\in(0, \pi/2)$
(resp. $x\in(\pi/2,\pi)$). There exists a unique $x_0$ in $(\pi/2, \pi)$ such that $h''(x_0)=0$.
Thus $h''(x)>0,\,x\in (0, x_0)$ and $h''(x)<0,\,x\in(x_0, \pi)$. This implies that $h'(x)$ is strictly increasing
in $x\in(0, x_0)$ and strictly decreasing in $x\in(x_0, \pi)$. As $h'(0)=0$, one has $h'(x_0)>0$.
By  $h'(\pi)= \pi^2-12 <0$, we get that there exists a unique $x_1$ in $(x_0, \pi)$ such that $h'(x_1)=0$.
We get that $h(x)$ is strictly increasing in $(0, x_1)$, and decreasing in $(x_1, \pi)$, with $h(0)=0$ and $h(\pi)=0$.
Thus $h(x)>h(0)=0,\, x\in(0, x_1)$ and $h(x)>h(\pi) =0,\,x\in(x_1, \pi)$. Hence, in all cases, one has  $h(x)>0$.
This completes the proof of the second inequality.

\vspace{.2cm}

It follows from the proof of the second inequality of Theorem \ref{thm0803} that $f<1/60$ and $f_1<0$. This implies the
following result.

\begin{corollary} For $x\in(0,\pi)$, we have
\begin{enumerate}
\item $M_1=\displaystyle\frac{2-x^4/60+\cos(x)}{3}<\frac{\sin(x)}{x}<\frac{8+7\cos(x)}{15-x^2}=M_2,$\\
\item $\displaystyle2\,\frac{x}{\sin(x)}+\frac{x}{\tan(x)}<3+\frac{x^4}{60}\frac{x}{\sin(x)}$,\\
\item $\displaystyle(8/7)\,\frac{x}{\sin(x)}+\frac{x}{\tan(x)}>\frac{15-x^2}{7}$.
\end{enumerate}
\end{corollary}
Recently, the following inequalities appeared in \cite{yang},
\begin{equation}\label{yanginequ}
Y_l=\cos\left(\frac{x}{\sqrt{x}}\right)^a<\frac{\sin(x)}{x}<\cos\left(\frac{x}{\sqrt{x}}\right)^{5/3}=Y_u
\end{equation}
for $x\in(0,\pi/2)$ with $a=\log(2/\pi)/\log(\cos(\sqrt{5}\pi/10))\approx 1.67141$. By using
Mathematica Software\textsuperscript{\tiny{\textregistered}} \cite{ru}, one can see
that
$$Y_l<M_1,\,x\in(0,1.06580)\quad {\rm and}\quad Y_u-M_2\in(0,-0.0009),\, x\in(0,\pi/2).$$

\vspace{.3cm}
\noindent{\bf Proof of Theorem \ref{wu-sri} \& \ref{1003}.} \rm
The proof of both theorems follow immediately from Theorem \ref{thm0803}.

\vspace{.3cm}
\noindent{\bf Proof of Theorem \ref{thm2}.} \rm
Let us consider the application
$$f(x)= \log\left(\frac{x}{\sin(x)}\right)-c \log (1/(1- x^2/\pi^2)),$$
for $x\in(0,\pi/2)$.
A simple computation gives
$$x\sin(x)(\pi^2-x^2)f'(x)= (\sin x-x\cos x)(\pi^2-x^2)-2cx^2\sin(x)=g(x).$$
One has
$$g'(x)/x= (\pi^2-2-4c)\sin( x)+(2-2c)x\cos(x)-x^2\sin(x) = h(x).$$
Finally,
$$h'(x)= (\pi^2-6c)\cos(x)-(4-2c)x\sin(x)-x^2\cos(x).$$
Now, if we select  $c=\pi^2/6$, then, as   $\pi^2-6c=0$ and $4-2c= 4-\pi^2/3>0$, we get $h'(x)<0$, so this finally leads to  $f(x)<0$. Hence, the first inequality follows.

For the proof of the second inequality, let $c=3/2$, then one has
$$g'(x)/x = (\pi^2-8)\sin (x)-x\cos(x) –(x^2)\sin(x)= h(x).$$
Since $h(\pi/2)<0,\, h(\pi/4)>0$, there exists an $x'$ in $(\pi/4, \pi/2)$ such that  $h(x')=0.$ We'll show that
$x'$ is unique. One has $h(x)/\sin (x)= \pi^2-8 -s(x)$, where $s(x)= x^2+ x/\tan(x)$.
Now, $$s'(x)\sin(x)^2 =  2x\sin(x)^2 + \cos (x)\sin( x)-x  =p(x).$$
Here $p'(x)= 4x\sin(x)\cos(x)>0$, which shows that $p(x)>p(0)=0$, so $s'(x)>0$, finally: $s(x)$ is a strictly increasing function. Thus the equation $s(x)= \pi^2-8$  has a single root (which is $x'$), so $h(x)>0$ for $x\in(0, x')$
and $h(x)<0$ for $x\in(x', \pi/2)$. Thus $g$ is increasing, resp. decreasing in the above intervals,
and $g(0)=g(\pi/2)=0$, so $g(x)>0$ for $x\in(0,\pi/2)$. This completes the proof of the second inequality.
$\hfill\square$

\begin{corollary}\label{cor1} For $x\in(0,\pi/2)$ one has
$$\frac{\pi ^2-x^2-\pi ^2 x^2/3}{\pi ^2-x^2}<\frac{x}{\tan(x)}<\frac{\pi ^2-4x^2}{\pi ^2-x^2}.$$
\end{corollary}

\begin{proof} After simplification the derivative of the function
$$f_c(x)= \log\left(\frac{x}{\sin(x)}\right)-c \log (1/(1- x^2/\pi^2)),$$
can be written as
$$f_c'(x)=\frac{1}{x}-{2c\,x}{\pi^2-x^2}-\frac{1}{\tan(x)}.$$
By the proof of Theorem \ref{thm2}, $f_{\pi^2/6}'(x)<0$ and $f_{3/2}'(x)>0$. Clearly,
$f_{\pi^2/6}=0=f_{3/2}$.
Now the proof of the inequalities is obvious.
\end{proof}
The second inequality in Corollary \ref{cor1} improves the first inequality in \eqref{zhu28c}.


\vspace{.3cm}
\noindent{\bf Proof of Theorem \ref{thm1}.} \rm For $t\in(0,1)$, let
$$f(t)=\frac{\pi  t \left(t^3+1\right)-\sin (\pi  t)}{\pi  t^3},\quad
g(t)=\frac{\pi t(2+t^4)-2\sin(\pi t)}{\pi t^3}.$$
By Theorem \ref{thm0803}, we get
\begin{eqnarray*}
f'(x)&=&\frac{3 \pi  t^3+\pi  \left(t^3+1\right)-\pi  \cos (\pi  t)}{\pi  t^3}-\frac{3 \left(\pi  t
   \left(t^3+1\right)-\sin (\pi  t)\right)}{\pi  t^4}\\
&=&\frac{3}{t^3}\left(\frac{\sin(\pi t)}{\pi t}-\frac{2-t^3+\cos(\pi t)}{3}\right)>0,
\end{eqnarray*}
and
\begin{eqnarray*}
g'(x)&=&2t-\frac{2}{t^2}\left(\frac{\cos(\pi t)}{ t}-\frac{\sin(\pi t)}{\pi t^2}\right)+
\frac{4}{t^3}\left(\frac{\sin(\pi t)}{ \pi t}-1\right)\\
&=&\frac{2}{t^3}\left(\frac{\sin(\pi t)}{ \pi t}-\frac{2-t^4+\cos(\pi t)}{3}\right)<0.
\end{eqnarray*}
Thus, the functions $f$ and $g$ are strictly increasing and decreasing in $t\in(0,1)$, respectively. Hence, the proof follows easily if we use the inequalities,
$$f(t)<\lim_{t\to 1}=2,\, g(t)>\lim_{t\to 1}g(t)=3,$$
and replace $t$ by $x/\pi$.$\hfill\square$

\vspace{.3cm}
\noindent{\bf Proof of Theorem \ref{gammathm}}. \rm  Utilizing \eqref{eulref}, the first inequality in Theorem \ref{thm1} is equivalent to
$$\frac{1+x/\pi-x/\pi}{\Gamma(1+x/\pi)\Gamma(1-x/\pi)}>\frac{\left(1-x/\pi\right)\left(1+x/\pi-
(x/\pi)^2\right)}{\Gamma(1+x/\pi-x/\pi)}.$$ Replacing $x$ by $x/\pi$
and $y$ by $1-x/\pi$ we get (1). Similarly, the second inequality in
Theorem \ref{thm1} can be written as,
$$\frac{1+x/\pi+1-x/\pi}{(1-(x/\pi)^2)(2+(x/\pi)^2)} < \frac{\Gamma(1+x/\pi)\Gamma(1-x/\pi)}{\Gamma(1+x/\pi-x/\pi)}.$$
If we replace $x$ by $1+x/\pi$ and $y$ by $1-x/\pi$ then we get the proof of part (2) for $1<x<2$. The rest of proof follows from Lemma \ref{lem1} .$\hfill\square$

%


\vspace{.3cm}
\noindent{\bf Proof of Theorem \ref{thm3}.} \rm Let $f(t)=f_1(t)/f_2(t),\,t\in(0,1/2)$, where
$$f_1(t)=1-3t^2-\cos(\pi t)\quad {\rm and}\quad f_2(t)=t^2(2+\cos(\pi t)).$$
A simple calculation gives
$$f'_1(t)=-6t-\pi\sin(\pi t)<0, \quad {\rm and}\quad f'_2(t)=t(4+2\cos(\pi t)-\pi t\sin(\pi t))>0.$$
Thus, $f$ is the product of two positive strictly decreasing functions,
this implies that $f$ is strictly decreasing in $t\in(0,1/2)$.
Applying l'H\^opital rule, we get
$\lim_{t\to 1/2}f(t)=1/2<f(t)<\alpha=\lim_{t\to 0}f(t)=\pi^2/6-1$. Here the first inequality implies
$$\frac{3(1-t^2)}{2+\cos(\pi t)}-1<\alpha t^2,$$
which is equivalent to
$$\frac{2+\cos(\pi t)}{3}>\frac{1-t^2}{1+\alpha t^2}.$$
Letting $\pi t=x\in(0,\pi/2)$, we get the second inequality of \eqref{thm3ineq},
and the third inequality of \eqref{thm3ineq} follows similarly. For the proof of the first inequality,
see \cite[(2.5)]{ael}. $\hfill\square$

\begin{corollary}\label{coro2} For $x\in(0,1)$, we have

\begin{equation}\label{stark}
\frac{4}{\pi}\frac{t}{1-t^2}<\tan\left(\frac{\pi t}{2}\right)<\frac{\pi}{2}\frac{t}{1-t^2}.
\end{equation}
\end{corollary}

\begin{proof} Let $f(t)=t/(\tan(\pi t/2)(1-t^2)),\,t\in(0,1/2)$. We get
$$f'(t)=\frac{(1+t^2)\sin(\pi t)-\pi t(1-t^2)}{2(1-t^2)^2\sin(\pi t/2)^2},$$
which is positive by \eqref{sinRed}. By l'H\^opital rule, we get
$\lim_{t\to 0}f(t)=2/\pi <f(t)<\pi/4=\lim_{t\to 1/2}f(t)$. This completes the proof.
\end{proof}
For $0<t<1$, letting $x=(\pi t)/2$ in \eqref{stark}, we get
\begin{equation}\label{starkineq}
\frac{8}{\pi^2-4x^2}<\frac{\tan(x)}{x}<\frac{\pi^2}{\pi^2-4x^2},
\end{equation}
which is so-called Becker-Stark inequality \cite{bs}. Thus, Corollary \ref{coro2} gives a simple proof of Becker-Stark inequality. It is easy to see that the second inequality in Corollary \ref{cor1} improves the first inequality
in \eqref{starkineq}.


\vspace{.3cm}
\noindent{\bf Proof of Theorem \ref{thm4}}. \rm Write $g(t)=g_1(t)/g_2(t),\ 0<t<1$, where
$$g_1(t)=\log\left(\frac{1+t^2}{1-t^2}\right),\,g_2(t)=\log\left(\frac{1}{\cos(\pi t/2)}\right).$$
We get,
$$\frac{g'_1(t)}{g'_2(t)}=\frac{8}{\pi}\frac{t}{(1-t^4)\tan(\pi t/2)}=\frac{8}{\pi}g_3(t).$$ One has,
$$g'_3(t)=\frac{\sin(\pi t)(1+3t^4)-\pi t(1-t^4)}{2\sin(\pi t/2)^2(1-t^4)^2},$$
which is negative by \eqref{zhu28a}.
Clearly $g_1(0)=0=g_2(0),$ thus $g$ is strictly decreasing by Lemma \ref{lem0avv}.
Using l'H\^opital rule, we get
$\lim_{t\to 0}g(t)=16/\pi^2 >g(t)>1=\lim_{t\to 1}g(t)$. Replacing $\pi t$ by $x$, we get the desired inequalities.$\hfill\square$

\begin{corollary}\label{coro3} For $x\in(0,\pi)$, one has
$$\left(\frac{\pi^2-x^2}{\pi^2+x^2}\right)^{4/(3\beta)}<\cos\left(\frac{x}{2}\right)^{4/3}<\frac{\sin(x)}{x}<
\left(\frac{\pi^2-x^2}{\pi^2+x^2}\right)^{4/(3\alpha)},$$
where $\alpha$ and $\beta$ are as in Theorem \ref{thm4}.
\end{corollary}

\begin{proof} The proof of the first inequality follows from Theorem \ref{thm4}, and the second inequality is
also well known \cite{neusana}. The third inequality is just \eqref{zhu28a}.
\end{proof}

\begin{theorem} For $x\in(0,\pi/2)$, we have
\begin{equation}\label{sinmax}
\frac{2}{\pi}+\frac{\pi-2}{\pi}\cos(x) < \frac{\sin(x)}{x},
\end{equation}
\begin{equation}\label{sinold}
\left(\frac{1+\cos}{2}\right)^{2/3} < \frac{\sin(x)}{x}< \frac{2\cdot2^{2/3}}{\pi}\left(\frac{1+\cos}{2}\right)^{2/3}
<\frac{4}{\pi}\frac{1+\cos(x)}{2}.
\end{equation}
\end{theorem}

\begin{proof}
The inequality \eqref{sinmax} may be rewritten as
$$f(x)= \pi\sin(x)-(\pi-2)x\cos(x)- 2x>0,$$
for $x\in(0,\pi/2)$.
Clearly, $f(0)=f(\pi/2)=0.$
On the other hand,

$$f(x)= 2(\cos (x)-1) +(\pi-2)x\sin(x) = -4\sin (x/2)^2 +2(\pi-2)x\sin (x/2)\cos (x/2).$$
 This implies that
$$f'(x)/(4\sin(x/2)\cos(x/2))= (\pi-2)/2-\tan( x/2) =g(x) $$
for $x>0$.
As $\tan (x/2)$ is strictly increasing, there is a unique $x'$ in $(0,\pi/2)$ such that $g(x')=0$. Also,
$g(x)>0$  for $x\in(0,x')$ and $g(x)<0$  for $x\in(x', \pi/2)$; i.e.  $x'$ is a maximum point of $f(x)$.
This gives $f(x)>f(0)=0$  for $x\in (0,x')$ and $f(x)>f(\pi/2)=0$  for $x \in (x',\pi/2)$. This completes the proof
of \eqref{sinmax}.

For the proof of \eqref{sinold}, let $$j(x)=\log\left(\frac{x}{\sin(x)}-\frac{2}{3}\frac{2}{1+\cos(x)}\right),$$
$x\in(0,\pi/2)$.
One has
\begin{eqnarray*}
j'(x)&=&\frac{1-x\cot(x)}{x}-\frac{2}{3}\frac{\sin(x)}{1+\cos(x)}\\
&=&\frac{1}{x}-\cot(x)-\frac{2}{3}\frac{1-\cos(x)}{\sin(x)}=j_1(x),
\end{eqnarray*}
which is negative, because the inequality  $j_1(x)<0$ can be written as $\sin(x)/x<(2+\cos(x))/3$,
which is so-called Cusa-Huygens inequality \cite{mit}. Thus,
$j(x)$ is strictly decreasing, and
$$\lim_{x->0}j(x)=0>j(x)>\log(\pi)-(5/3)\log(2)=\lim_{x\to \pi/2}j(\pi/2)\approx -0.010515.$$
By these inequalities we get
$$\cos(x/2)^{4/3}<\frac{\sin(x)}{x}<\exp((5/3)\log(2)-\log(\pi))\cos(x/2)^{4/3}.$$
The proof of the first and second inequality is completed, and the proof of the third inequality is trivial.
\end{proof}

The inequality \eqref{sinmax} improves the following one
\begin{equation}\label{sandor1702}
1-2\frac{\pi-2}{\pi^2}x<\frac{\sin(x)}{x},\quad 0<x<\pi,
\end{equation}
which was proved in \cite{sandor17} as an application of the concavity
of $\sin(x)/x$. Indeed, the inequality
$$1-2\frac{\pi-2}{\pi^2}x<\frac{2}{\pi}+\frac{\pi-2}{\pi}\cos(x),$$
is equivalent to
\begin{equation}\label{kober1702}
\cos(x)>1-\frac{2}{\pi}x,\quad 0<x<\frac{\pi}{2},
\end{equation}
which is Kober's inequality, see \cite{mit,sandor17b}.



%

\vspace{.3cm}
\noindent{\bf Proof of Theorem \ref{thm1702}.} \rm For $x\in(0,\pi/2)$, let
$$f(x)=\frac{x}{\sin(x)}\frac{x(\pi-x)}{(2x-\pi)^2}-\frac{4x^2}{(2x-\pi)^2}.$$
We get,
$$f'(x)=-\frac{16\pi^2\cos(x)}{(2x-\pi)^3\sin(x)^2}g(x),$$
where
$$g(x)=\tan(x)-\frac{\sin(x)^2}{\cos(x)}-x.$$
One has,
$$g'(x)=\frac{\sin(x)^2}{\cos(x)}(\sin(x)^2-\sin(x)-2)<0,$$
as $0<\sin(x)<1$, and $g(x)=0$. Thus, $g<0$, and in result $f'(x)<0$, this implies that
$f$ is strictly decreasing. By l'H\^opital rule we get
$$\lim_{x\to \pi/2}f(x)=\frac{\pi^2-8}{2}<f(x)<\lim_{x\to 0}f(x)=\frac{16}{\pi}-4,$$
this implies the proof. $\hfill\square$

\vspace{.5cm}
Replacing $x$ by $\pi/2-x$ in the inequalities of Theorem \ref{thm1702}, we get the following corollary as a result.

\begin{corollary}\label{cor1702} For $x\in(0,\pi/2)$, we have
$$\frac{16(\pi-x)}{4\pi^2+\alpha (2x-\pi)^2}<\frac{\sin(x)}{x}<\frac{16(\pi-x)}{4\pi^2+\beta (2x-\pi)^2},$$
where $\alpha$ and $\beta$ are as in Theorem \ref{thm1702}.
\end{corollary}


\end{document}



%

\maketitle

\begin{abstract}
In this paper we study the two sided inequalities of the trigonometric and hyperbolic functions.
\end{abstract}

%


\section{Introduction}

Since last one decate many authors got interest to study the following inequalities

\begin{equation}\label{laztri}
(\cos x)^{1/3}<\frac{\sin x}{x}<\frac{\cos x+2}{3},
\end{equation}
for $0<|x|<\pi/2$. The first inequality is due to by D.D. Adamovi\'c and D.S. Mitrinovi´c \cite[p.238]{mit}, and the second one was obtained
by N. Cusa and C. Huygens \cite{sanben}. The hyperbolic version of \eqref{laztri} is given as 
\begin{equation}\label{lazhyp}
(\cosh x)^{1/3}<\frac{\sinh x}{x}<\frac{\cosh x+2}{3},
\end{equation}
for $x\neq 0$. The first inequality in \eqref{lazhyp} was obtained by Lazarevi\'c \cite[p. 270]{mit}, 
and the second inequality is called sometime
hyperbolic Cusa-Huygens inequality \cite{neusan}.

There are many results on the refinement of the inequalities \eqref{laztri} and \eqref{lazhyp} in the literature,
 e.g, see 
\cite{avv2,neu,neu1,neu2,neusan,nejoz1,joz1,joz111,joz2,joz3,sanben} and the references therein. 
In this paper, we give the simple proofs of some known results, and establish the new inequalities as well as.
The main result of this paper reads as follows.

\begin{theorem}\label{thm0}
For $x\in(0,\pi/2)$, the following inequalities 
\begin{equation}\label{tanx2}
\frac{x^{3/2}}{2(\sin x)^{1/2}}<\tan \frac{x}{2}<\frac{x^{2/\alpha}}{2(\sin x)^{2/\alpha-1}}
\end{equation}
hold true, where $\alpha=2\log(\pi/2)/\log(2)\approx 1.30299.$
\end{theorem}

\begin{theorem}\label{thmm} For $x\in(0,\pi/2)$, the following inequalities
$$
\frac{x}{\tan (x/2)}-1<\exp\left(\frac{1}{2}\left(\frac{x}{\tan x}-1\right)\right)
< \frac{\sin x}{x}< 
\exp\left(\left(\log\frac{\pi}{2}\right)\left(\frac{x}{\tan x}-1\right)\right)$$
hold true.
\end{theorem}

\begin{theorem}\label{thmm1}
For $x\in(0,\infty)$, the following inequalities 
$$\frac{x}{\tanh (x/2)}-1<\exp\left(\frac{1}{2}\left(\frac{x}{\tanh x}-1\right)\right)
< \frac{\sinh x}{x}< 
\exp\left(\left(\frac{x}{\tanh x}-1\right)\right)$$
hold true.
\end{theorem}
The second inequality in Theorem \ref{thmm} and \ref{thmm1} is known \cite{joz1}. For these, however,  here a new method of proof is offered.
\section{Preliminaries}

In this section we give few lemmas, they will be used in the proof of theorems.

For $|x|<\pi$, the following power series expansions can be found in \cite[1.3.1.4 (2)-(3)]{jef},
\begin{equation}\label{xcot}
x \cot x=1-\sum_{n=1}^\infty\frac{2^{2n}}{(2n)!}|B_{2n}|x^{2n},
\end{equation}

\begin{equation}\label{cot}
\cot x=\frac{1}{x}-\sum_{n=1}^\infty\frac{2^{2n}}{(2n)!}|B_{2n}|x^{2n-1},
\end{equation}
and 
\begin{equation}\label{coth}
{\rm \coth x}=\frac{1}{x}+\sum_{n=1}^\infty\frac{2^{2n}}{(2n)!}|B_{2n}|x^{2n-1},
\end{equation}
where $B_{2n}$ are the even-indexed Bernoulli numbers, see \cite[p. 231]{IR}. 
The following expansions can be obtained directly from (\ref{cot}) and (\ref{coth}),

\begin{equation}\label{cosec}
\frac{1}{(\sin x)^2}=-(\cot x)'=\frac{1}{x^2}+\sum_{n=1}^\infty\frac{2^{2n}}{(2n)!}
|B_{2n}|(2n-1)x^{2n-2},
\end{equation}

\begin{equation}\label{cosech}
\frac{1}{(\sinh x)^2}=-({\rm coth} x)'=\frac{1}{x^2}-\sum_{n=1}^\infty\frac{2^{2n}}{(2n)!}(2n-1)|B_{2n}|x^{2n-2}.
\end{equation}

The following result is due to Biernacki and Krzy\.z \cite{bk}, which will
be very useful in studying the monotonicity of certain power series. This result was also simply proved in
\cite{hvv} by Heikkala et al.

\begin{lemma}\label{lembk} For $0<R\leq \infty$. Let $A(x)=\sum_{n=0}^\infty a_nx^n$ and 
$B(x)=\sum_{n=0}^\infty b_nx^n$ be two real power series converging on the interval $(-R,R)$. If the sequence
$\{a_n/b_n\}$ is increasing (decreasing) and $b_n>0$ for all $n$, then the function $A(x)/B(x)$ is also
increasing (decreasing) on $(0,R)$.
\end{lemma}


In \cite{avv1}, Anderson et al. proved the following result, which is sometime called Monotone l'H\^opital rule.
\begin{lemma}\label{lem0}
For $-\infty<a<b<\infty$,
let $f,g:[a,b]\to \mathbb{R}$
be continuous on $[a,b]$, and be differentiable on
$(a,b)$. Let $g^{'}(x)\neq 0$
on $(a,b)$. If $f^{'}(x)/g^{'}(x)$ is increasing
(decreasing) on $(a,b)$, then so are
$$\frac{f(x)-f(a)}{g(x)-g(a)}\quad and \quad \frac{f(x)-f(b)}{g(x)-g(b)}.$$
If $f^{'}(x)/g^{'}(x)$ is strictly monotone,
then the monotonicity in the conclusion
is also strict.
\end{lemma}

\begin{lemma}\label{lem1a} (1) The function
$$f(x)=x+\sin x-\frac{x^2}{\tan (x/2)}$$
is strictly increasing from $(0,\pi/2)$ onto $(0,c),\,c=(4+2\pi-\pi^2)/4\approx 0.10336$. In paticular,
$$1+\frac{\sin x}{x}-\frac{c}{x}<\frac{x}{\tan (x/2)}< 1+\frac{\sin x}{x},$$
for $x\in(0,\pi/2)$.\\
(2) The function 
$$g(x)=\frac{1}{x}\left(x\frac{\cos x}{\sin x}-1\right)\left(\frac{\cos x}{\sin x}+\frac{1}{\sin x}\right)$$
is strictly increasing from $(0,\pi/2)$ onto $(-4/3,-4/\pi)$.
\end{lemma}

\begin{proof} For the proof of (1), we get 
$$f'(x)=\frac{1}{2}\left(\frac{1}{\sin (x/2)}\right)^2(x-\sin x)^2>0,$$
the limiting values are clear.
For (2), derivation gives
$$g'(x)=\frac{1}{x^2(1-\cos x)}\left(x+\sin x-\frac{x^2}{\tan (x/2)}\right),$$
which is positive by part (1), and limiting values follow from l'H\^opital rule. This completes the proof.
\end{proof}

\begin{lemma}\label{lem3} (1) The function
$$f(x)=\cos x-\left(1-\frac{x^2}{3}\right)^{3/2}$$
is strictly decreasing from $(0,\pi/2)$ onto $(0,c_1)$, where $c_1=-(12-\pi^2)^{3/2}/(24\sqrt{3})\approx -0.07480$. In particular, for $x\in(0,\pi/2)$
\begin{equation}\label{cos}
\left(1-\frac{x^2}{3}\right)^{3/2}-c_1<\cos x < \left(1-\frac{x^2}{3}\right)^{3/2},
\end{equation}
(2) For $x\in(0,\pi/2)$, we have
\begin{equation}\label{tanx1}
\frac{x}{4}\left(\frac{5}{2}-\frac{x}{2\tan x}\right) < \tan\frac{x}{2}.
\end{equation}
(3) For $x\in(0,\infty)$, we have
\begin{equation}\label{tanhx1}
\frac{x}{4}\left(\frac{5}{2}-\frac{x}{2\tanh x}\right) < \tanh\frac{x}{2}.
\end{equation}
\end{lemma}

\begin{proof} By utilizing the inequality $(1-x^2/6)<(\sin x)/x$ from \cite[Theorem 3.1]{kvv}, we get
\begin{eqnarray*}
f'(x)&=&x(1-x^2/3)^{1/2}-\sin x\\
&<&x\left((1-x^2/3)^{1/2}-(1-x^2/6)\right)=x\,g(x),
\end{eqnarray*}
which is negative, because
$$g'(x)=\frac{x}{3}\left(1-\frac{1}{(1-x^2/3)^{1/2}}\right)<0,$$
and $g(x)=0$. This implies that $f$ is decreasing, and limiting values are clear.

For the proof of part (2), we use the indentity $\tan (x/2)=(1-\cos x)/\sin x$ and write the inequality \eqref{tanx1} 
as 
\begin{equation}\label{newsin}
\frac{\sin x}{x} < \frac{8(1-\cos x)+x^2\cos x}{5 x^2},\quad x\in (0,\pi/2).
\end{equation}
Let 
$$h(x)=5 x\sin x-8(1-\cos x)+x^2\cos,\quad x\in(0,\pi/2).$$
Applying the inequality $(\cos x)^{1/3}< \sin x/x,$ we get
\begin{eqnarray*}
h'(x)&=&3x\cos x-(3-x^2)\sin x< x(3\cos x-(3-x^2)(\cos x)^{1/3})\\
&=& 3x (\cos x)^{1/3}\left((\cos x)^{2/3}-\left(1-\frac{x^2}{3}\right)\right),
\end{eqnarray*}
which is negative by part (1), and $h(0)=0$. The proof of part (3) is similar to the proof of part (2), if we use
the indentity 
$\tanh x/2=(\cosh x-1)/\sinh x$.
\end{proof}


\section{Proofs of the theorems}

\subsec{\bf Proof of Theorem \ref{thm0}.} \rm Let $f(x)=f_1(x)/f_2(x)$,  where
 $$f_1(x)=\log(x/\sin x),\, f_2(x)=\log(1/\sqrt{(1+\cos x)/2}),$$
and clearly $f_1(x)=f_2(x)=0$. We get
$$\frac{f_1'(x)}{f_2'(x)}=-\frac{1}{x}\left(x\frac{\cos x}{\sin x}-1\right)\left(\frac{\cos x}{\sin x}+\frac{1}{\sin x}\right),$$ 
which is strictly decreasing by Lemma \ref{lem1a}(2). Hence the function $f$ is strictly decreasing by Lemma \ref{lem0},
and by l'H\^opital rule $\lim_{x\to 0}f(x)=4/3$ and $\lim_{x\to \pi/2}f(x)=2\log(\pi/2)/\log(2)=1.30299\ldots=\alpha$. 
This implies the following inequalities
\begin{equation}\label{sinx}
\left(\sqrt{\frac{1+\cos x}{2}}\right)^{4/3}<\frac{\sin x}{x}<\left(\sqrt{\frac{1+\cos x}{2}}\right)^\alpha.
\end{equation}
By utilizing the indentity $\tan x/2=(\sin x)/(1+\cos x)$ we get \eqref{tanx2}.$\hfill\square$

\vspace{.4cm}
In \cite{joz11}, it is proved that
\begin{equation}\label{3001}
\frac{\sin x}{x}<\sqrt{\frac{1+\cos x}{2}}\frac{\sqrt{(1+\cos x)/2}+2}{3}=\frac{t^{1/2}(t^{1/2}+2)}{3}
=\frac{t+2t^{1/2}}{3},
\end{equation}
if one denotes $t=(\cos x+1)/2$.
The right side of \eqref{sinx} cannot be
compared with inequality \eqref{3001}, and
this inequality improves the Cusa-Huygens inequality $(\sin x)/x< (\cos
x+2)/3 = (2t+1)/3$.

The following theorem is known \cite{joz2}. Here another proof appears, but both proofs use Lemma \ref{lem0}.

\begin{theorem}\label{th2}
The function 
$$f(x)=\frac{\log((\sinh x)/x)}{\log(\cosh (x/2))}$$
is strictly decreasing from $(0,\infty)$ onto $(4/3,2)$. In particular, 
$$\left(\frac{1+\cosh x}{2}\right)^{2/3}<\frac{\sinh x}{x}<\left(\frac{1+\cosh x}{2}\right).$$
\end{theorem}

\begin{proof}  Let $f_1(x)=\log((\sinh x)/x),\,f_2(x)=\log((\cosh x)/x)$, and $f_1(0)=f_2(0)=0$. We get
$$\frac{f_1'(x)}{f_2'(x)}=2g_1(x),$$
where $$g_1(x)=x\left(1+\frac{\cosh x}{\sinh x}\right)\frac{1}{(\sinh x)^2}
\left(\frac{\cosh x}{x}+\frac{\sinh x}{x}\right).$$ Differentiation gives,
$$g_1'(x)=\frac{1}{2x^2}g_2(x),$$ where 
$$g_2(x)=2\frac{\cosh (x/2)}{\sinh (x/2)}+\frac{x}{(\sinh (x/2))^2}\left(1-x\frac{\cosh (x/2)}{\sinh (x/2)}\right),$$
which is positive because
$$g_2'(x)=\frac{x}{2(\sinh (x/2))^4}(x(2+\cosh x)-3\sinh x)>0$$
by inequalitiy $(\sinh x)/x<(2+\cosh x)/3$, and $g_2$ tends to $0$ as $x$ tends to $0$. This implies that 
$f_1'/f_2'$ is increasing, Hence by Lemma \ref{lem0} $f$ is increasing, and limiting values follow 
from l'H\^opital rule. This completes the proof.
\end{proof}

\begin{corollary}\label{cor1} For $x\in(0,\infty)$, we have
$$(\cosh x)^{1/3}<\left(\frac{1+\cosh x}{2}\right)^{2/3}<\frac{\sinh x}{x}<\frac{1+\cosh x}{2}
<\left(\frac{2+\cosh (x/2)}{3}\right)^4.$$
\end{corollary}

\begin{proof} The proof of the first inequality is trivial, because it can be simplified as $0<(1-\cosh x)^2$. 
The second and third inequality follow from Theorem \ref{th2}. For the proof of the fourth inequality, it is enough to prove that
the function 
$$h(y)=3\left(\frac{y+2}{3}\right)^{1/4}-\left(\frac{y+1}{2}\right)^{1/2}-2$$
is negative for $y\in(0,\infty)$. We get
$$h'(x)=\frac{3^{3/4}}{4(1+y)^{1/2}}\left(h_1(y)-\frac{\sqrt{2}}{3^{3/4}}\right),$$
where
$$h_1(x)=\frac{(1+y)^{1/2}}{(2+y)^{3/4}}.$$ Now we prove that $h_1$ is negative, we get
$$(\log(h_2(y)))'=-\frac{y-1}{4(2+3y+y^2)}<0,$$
and $h_2(1)=\sqrt{2}/3^{3/4}$. This implies that $h$ is decreasing, and $h$ tends to $0$ when $y$ tends to 1, 
hence the proof of fourth inequality is completed.
\end{proof}

\subsec{\bf Proof of Theorem \ref{thmm}.} \rm  It is well known that  $L(a,b)< (a+b)/2$, where $L(a,b)=(b-a)/(\log b-\log a)$ denotes the Logarithmic mean of two distint real numbers $a$ and $b$. Particularly, $L(t,1)<(t+1)/2$, for $t>1$. This inequality can be written as
$\log t>2(t-1)/(t+1)$. Letting $t=\tan(x/2)/ (x-\tan
(x/2))>1$ and using \eqref{tanx1}, we get
$$\log\left(\frac{\tan (x/2)}{x-\tan (x/2)}\right)>\frac{2(\tan (x/2)-x)}{x}>\frac{\tan x-x}{2\tan x},$$
this implies the proof of first inequality.

For the proof of second and third ienquality, we define $f(x)=f_1(x)/f_2(x),\,x\in(0,\pi/2)$, where
$f_1(x)=\log(x/\sin x),\ f_2(x)=1-x/\tan x,$
and clearly $f_1(0)=f_2(0)=0$.
Differentiation with respect $x$ gives
$$\frac{f_1'(x)}{f_2'(x)}=\frac{A(x)}{B(x)},$$
where $$A(x)=1-x \frac{\cos x}{\sin x}\quad {\rm and }\quad B(x)=\frac{x^2}{\sin x^2}-x \frac{\cos x}{\sin x}.$$
By \eqref{xcot} and \eqref{cosec}, we get 
$$A(x)=\sum_{n=1}^\infty\frac{2^{2n}}{(2n)!}|B_{2n}|x^{2n}=\sum_{n=1}^\infty a_nx^{2n},$$
and 
$$B(x)=\sum_{n=1}^\infty\frac{2^{2n}2n}{(2n)!}|B_{2n}|x^{2n}=\sum_{n=1}^\infty b_nx^{2n}.$$
The function $c_n=a_n/b_n=1/(2n)$ is decreasing in $n=1,2,\ldots$. This implies from Lemma \ref{lembk} that $f_1'/f_2'$ is decreasing, and the function $f$ is decreasing in $x\in(0,\pi/2)$ by Lemma \ref{lem0}.
By l'H\^opital rule, we get $\lim_{x\to 0}f(x)=1/2$ and $\lim_{x\to \pi/2}f(x)=\log(\pi/2)$, this completes the proof of the second and the third inequality.$\hfill\square$

\subsec{\bf Proof of Theorem \ref{thmm1}.} \rm The proof of the first inequality is similar to the 
proof of the first inequality of Theorem \ref{thmm} if we use 
(\ref{tanhx1}). For the proof of second and third inequality, we define the function
$g(x)=g_1(x)/g_2(x),\,x\in(0,\infty)$, where 
$g_1(x)=\log(\sinh(x)/x)$ and $g_2(x)=x/\tanh(x)-1$, and $g_1(0)=g_2(0)=0$. Differentiating $g$ with respect $x$
we get
$$g'(x)=\frac{A(x)}{B(x)},$$
where
$A(x)=1-x\,{\rm coth}(x)$ and $B(x)=(x/\sinh x)^2-x\,{\rm coth}x$. Now the monotonicty of $g$
 follows from (\ref{cot}), (\ref{cosech}), Lemma
\ref{lembk} and \ref{lem0}. 
Applying the l'H\^opital rule, we get $\lim_{x\to 0}g(x)=1/2$ and $\lim_{x\to \infty}=1$, this completes the proof.
$\hfill\square$

The second inequality of Theorem \ref{thmm} can be written as
$$\log\left(\frac{x}{\sin x}\right)<\frac{\sin x-x\cos x}{2\sin x},$$
and this appears in \cite{joz1}.
 The similar is true for the second inequality of Theorem \ref{thmm1}, which appears in \cite{joz1}
as 
$$\log\left(\frac{\sinh x}{x}\right)<\frac{x\cosh x -\sinh x}{2\sinh x}.$$

\begin{corollary}\label{corjoz} The following relations hold true,
\begin{enumerate}
\item $\displaystyle x\tan\left(\frac{x}{2}\right)<\log\left(\frac{1}{\cos x}\right)
<\displaystyle\frac{x}{2}\tan x,\quad x\in(0,\pi/2),$
\item $\displaystyle\frac{x}{2}\tanh x<\log(\cosh x)<\displaystyle\frac{x^2}{2},\quad x\in(0,\infty),$
\item $\displaystyle\frac{x}{2}\tanh x<\log(\cosh x)<\displaystyle\frac{\sinh x \tanh x}{2},\quad x\in(0,\infty).$
\end{enumerate}
\end{corollary}

\begin{proof} 
Applying the Hadamard inequalities (see \cite{joz111}) to the convex
function $\tan t$, we get:
$$
x\tan\left(\frac{x+0}{2}\right)=x\tan (x/2)<\int_0^x\tan t\, dt<\frac{x(\tan x+\tan x)}{2}=x\tan (x/2),
$$
so (1) follows.
For part (2), it is sufficient to prove that
the functions  $m(x)= \log (\cosh x)- (x/2)\tanh x $ and $n(x)= x^2/2 – \log (\cosh x) $  are strictly increasing for $x>0$.
This follows from $m'(x)=(\sin x-x/\cos x)/(2\cos x)>0$ and $n'(x)=1-1/(\cosh x)^2>0$.
It is interesting to note that, the right side of (2) cannot be compared
to right side of (3).
\end{proof}

%
%
%
%
%


\begin{lemma}\label{jozlem1} The following inequalities hold,
\begin{enumerate}
\item $\displaystyle\log\left(\frac{1}{\cos x}\right)<\displaystyle\frac{1}{2}\sin x \tan x,
\quad x\in\left(0,\frac{\pi}{2}\right),$
\item $\displaystyle\log\left(\frac{x}{\sin x}\right)>\displaystyle\frac{\sin x-x\cos x}{2x},
\quad x\in\left(0,\frac{\pi}{2}\right),$
\item $\displaystyle\log\left(\frac{\sinh x}{x}\right)
<\displaystyle\frac{x\cosh x-\sinh x}{2x},\quad x>0,$
\item $\displaystyle\log\left(\frac{\sinh x}{x}\right)>\displaystyle\frac{x\cosh x-\sinh x}{2\sinh x},\quad x>0.$
\end{enumerate}
\end{lemma}

\begin{proof} Letting $$t(x)=\frac{\sin x\tan x}{2}-\log\left(\frac{1}{\cos x}\right),$$
we get $t'(x)=\sin x(\cos x-1)^2/(\cos x)^2>0$. Hence $t$ is strictly increasing, this implies (1). 
As the logarithmic mean of $t$ and 1 is
$L(t,1)=(t-1)/\log t$, from $\log t>1-1/t$ for $t>1$, we get $L(t,1)<t$. Putting $t=x/\sin x$, we obtain
$\log (x/\sin x)>1-(\sin x)/x$. Thus, it is sufficient to prove that 
$$1-\frac{\sin x}{x}>\frac{\sin x}{2x}-\frac{\cos x }{2},$$
which is equivalent to write $(\sin x)/x<(\cos x+2)/3$. This is a known inequality, namely the so-called Cusa-Huygens inequality \eqref{laztri}, thus (2) follows. For the proof of (3) we use the similar method as in the circular case:
By inequality  $\log t< t-1$ for $t>1$ applied to $t= (\sinh x)/x$, to deduce (3) it
is sufficient  to prove that:
$$\frac{\sinh x}{x} -1< \frac{\cosh x}{2}- \frac{\sinh x}{2x}.$$
 But this is equivalent to
$(\sinh x)/x  < (\cosh x+2)/3$,
which is known as the hyperbolic Cusa-Huygens inequality \eqref{lazhyp}. For the proof of (4), see \cite{joz2}.
\end{proof}

The last inequality of Theorem \ref{thmm1} can be written as
\begin{equation}\label{nnn}
\displaystyle\log\left(\frac{\sinh x}{x}\right)<\displaystyle\frac{x\cosh x-\sinh x}{\sinh x}.
\end{equation}
We observe that the inequality (3) in Lemma \ref{jozlem1} can not be compared with \eqref{nnn}.

\begin{theorem} The following functions 
$$F_1(x)= \frac{\log (x/\sin x)}{\log (\cosh x)},\quad
F_2(x)= \frac{\log (\cosh x)}{\log (\cos x)},\quad
F_3(x)= \frac{\log (x/\sin x)}{\log (\cos x)}$$  
are strictly increasing in $x\in(0,\pi/2)$.
\end{theorem}
The proof of the functions $F_1$ can be found in \cite{nejoz1}. 
The proof of $F_2$ is a consequence of Corollary \ref{corjoz}, by considering the derivative  
$$F'_2(x)=  \tanh x \log(\cos x)+ \cos x \log(\cosh x) >0.$$ The proof for the $F_3$
can be found in \cite{neu1}, here we give a simple proof. 
Let $F_3(x)=u(x)/v(x)$, the proof follows if $(u(x)/v(x))'$ is negative or one has to show that
$u'(x)v(x)<u(x)v'(x)$. This is equivalent to 
\begin{equation}\label{F3}
\log\left (\frac{1}{\cos x}\right)\left(\frac{\sin x-x \cos x}{x \sin x}\right)<\log\left(\frac{x}{\sin x}\right)\tan x .
\end{equation}
The inequalitiy \eqref{F3} follows immediately from Lemma \eqref{jozlem1} (1) and (2). This completes the proof.

\begin{corollary}\label{jozcor} For $x\in(0,\pi/2)$, one has
$$\frac{3\sin x}{2+\cos x}<\frac{8\sin(x/2)-\sin x}{3}<x<\frac{\sin x}{((\cos x+1)/2)^{2/3}}.$$
\end{corollary}

\begin{proof}
By applying the Mitrinovic-Adamovic inequality $(\sin t)/t > (\cos t)^{1/3}$ for
$t=x/2$, we get
$$\left(\frac{2\sin (x/2)\cos (x/2)}{x}\right)^3  > (\cos x/2)^4 $$ 
by multiplying both sides with
$\cos (x/2)$.
Since $\cos (x/2)= (1+ \cos x)^{1/2}$,
we get the inequality:
$$( \sin x)/x>  ((\cos x+1)/2)^{2/3}.$$
By putting $x/2$
in the Cusa-Huygens inequality  $(\sin x)/x< (\cos x+2)/3$, we get
\begin{equation}\label{josre}
\frac{\sin (x/2)}{x/2}< \frac{\cos (x/2) +2}{3}.
\end{equation}
Clearly, $(\cos x-1)^2>0,\,t=\cos (x/2)$, this implies
\begin{equation}\label{josre1}
\frac{\cos t+2}{3}< \frac{3}{4-\cos t}.
\end{equation}
By \eqref{josre} and \eqref{josre1}, one has 
$$\sin (x/2)<\frac{3x}{2(4-\cos (x/2))},$$
this is  stronger than the Cusa-Huygens inequality  
$(\sin x)/x<(2+\cos x)/3$.
Indeed, we will prove that:
\begin{equation}\label{josre2}
\frac{3\sin x}{2+\cos x}<\frac{1}{3} (8\sin (x/2) -\sin x).
\end{equation}
By letting $t=\cos (x/2)$ in \eqref{josre2}, after some simple transformations, we have to
prove that
$2t^3-8t^2+10t-4<0$, this is equivalent to $(t-1)^2(t-2)<0$, which is obvious.
\end{proof}

\begin{remark}\rm
The first two inequalities of Corollary \ref{jozcor} offer an improvement of the
Cusa-Huygens inequality, while the last inequality is in fact the first
relation of (\ref{sinx}), with an elementary proof.
The authors have found in a book on history of number $\pi$,  that the
second inequality of Corollary \ref{jozcor} was discovered by geometric and
heuristic arguments (see \cite{campan}) by C. Huygens.
\end{remark}



\font\fFt=eusm10 
\font\fFa=eusm7  
\font\fFp=eusm5  
\def\K{\mathchoice
{\hbox{\,\fFt K}}
{\hbox{\,\fFt K}}
{\hbox{\,\fFa K}}
{\hbox{\,\fFp K}}}
\def\T{\mathchoice
{\hbox{\,\fFt T}}
{\hbox{\,\fFt T}}
{\hbox{\,\fFa T}}
{\hbox{\,\fFp T}}}

%

\begin{abstract}
In this paper authors refine the Carlson's inequalities for inverse cosine function,
and the Shafer's inequalities for inverse tangent function.
\end{abstract}
\maketitle

%

\section{Introduction}

During the past fifteen years, numerous authors have studied various inequalities for trigonometric
functions \cite{avv1,kvv,mal1,mal2,neusan,nejoz1,joz1,joz2,joz3,sandor17a,sz,zredh}. Thus, some classical and also more recent inequalities, such as inequalities of Jordan, Cusa-Huygens, Shafer-Fink, and Wilker have been refined and generalized.
One of the key methods
in these studies has been so called monotone l'Hospital Rule from \cite{avv1} and an extensive survey of the related
literature is given in \cite{avz}. This Rule is formulated as Lemma \ref{lem0}  and it will be also applied here.
Motivated by these studies,
in this paper we make a contribution to the topic by sharpening the Carlson's and Shafer's inequalities, and our inequalities refine the existing results in the literature.

We start our discussion with the following well-known inequalities,
\begin{equation}\label{cusa}
\cos(t)^{1/3}<\frac{\sin(t)}{t}<\frac{\cos(t)+2}{3},
\end{equation}
for $0<|t|<\pi/2$. The first inequality is due to by D.D. Adamovi\'c
and D.S. Mitrinovi\'c \cite[p.238]{mitr}, while the second inequality was
obtained by N. Cusa and C. Huygens \cite{huy}. These inequalities
can be written as,
$$
\frac{3\sin(t)}{2+\cos(t)}<t<\frac{\sin(t)}{\cos(t)^{1/3}}.
$$
For the further studies and refinements of inequalities in
\eqref{cusa}, e.g., see \cite{avz,kvv,neusan,nejoz1,joz1,joz2,gendi}
and the references therein. For the easy references we recall the
following inequalities
\begin{equation}\label{refine}
\cos\left(\frac{t}{2}\right)^{4/3}<\frac{\sin(t)}{t}<\cos\left(\frac{t}{3}\right)^3,
\end{equation}
the first inequality holds for $t\in(0,\pi/2)$ \cite{nejoz1},
while the second one is valid for
$t\in(-\sqrt{27/5},\sqrt{27/5})$, and was proved by Kl\'en et al.
\cite{kvv}. The first inequality in \eqref{refine} refines the following one
$$t<\frac{2\sin(t)}{1+\cos(t)},\quad 0<x<\frac{\pi}{2},$$
which was constructed in \cite{qicu} by using the Chebyshev's
integral inequality.


The Oppenheim's problem \cite{ogi,oppe,mitr} states that: to determine the greatest $a_2$ and least $a_3$ 
as a function of $a_1>0$, such that the following inequalities
\begin{equation}\label{oppeineq}
\frac{a_2\sin(x)}{1+a_1\cos(x)}<x<\frac{a_3\sin(x)}{1+a_1\cos(x)},
\end{equation}
hold for all $x\in(0,\pi/2)$.
A partial solution of this problem was given by Oppenheim and Carver \cite{oppe}, they showed that
\eqref{oppeineq} holds for all $a_1\in(0,1/2)$ and $x\in(0,\pi/2)$ when $a_2=1+a_1$ and $a_3=\pi/2$.
 In 2007, Zhu \cite[Theorem 7]{z4} solved the Oppenheim's problem completely by proving that the inequalities in 
\eqref{oppeineq}
hold if $a_1,\,a_2$ and $a_3$ are as follows:
\begin{enumerate}
\item if $a_1\in(0,1/2)$, then $a_2=1+a_1$ and $a_3=\pi/2$,\\
\item if $a_1\in(1/2,\pi/2-1)$, then $a_2=4a_1(1-a_1^2)$ and $a_3=\pi/2$,\\
\item if $a_1\in(\pi/2-1,2/\pi)$, then $a_2=4a_1(1-a_1^2)$ and $a_3=1+a_1$,\\
\item if $a_1>2/\pi$, then $a_2=\pi/2$ and $a_3=1+a_1$,\\
\end{enumerate}
where $a_2$ and $a_3$ are the best possible constants in (1) and (4), while $a_3$
is the best possible constant in (2) and (3). Thereafter, the Carver's solution was
extended to the Bessel functions for the further results by Baricz \cite{bar1,bar2}.
On the basis of computer experiments we came up that the following lower and upper bounds for $x$,
\begin{equation}\label{zhuineq}
\frac{(\pi/2)\sin(x)}{1+(2/\pi)\cos(x)}<x<\frac{\pi\sin(x)}{2+(\pi-2)\cos(x)}
\end{equation}
are the best possible bounds, and can be obtained from case (4) and (3), respectively.

Recently, Qi et al. \cite{qiluo} have given a new proof of Oppenheim's problem, and deduced the following inequalities,
\begin{equation}\label{qioppeineq}
\frac{(\pi/2)\sin(x)}{1+(2/\pi)\cos(x)}<x<\frac{(\pi+2)\sin(x)}{\pi+2\cos(x)}.
\end{equation}
for $x\in(0,\pi/2)$.
It is obvious that
$$((\pi-2)-4)(1-\cos(x))<0,$$
which is equivalent to
$$\frac{\pi\sin(x)}{2+(\pi-2)\cos(x)}<\frac{(\pi+2)\sin(x)}{\pi+2\cos(x)}.$$
This implies that the second inequality of \eqref{zhuineq} is better than the corresponding inequality of
\eqref{qioppeineq}.

Our first main result reads as follows, which refines the inequalities in \eqref{zhuineq}.

\begin{theorem}\label{thm1} For $x\in(0,\pi/2)$, we have
\begin{equation}\label{thm1inequ}
C_\alpha< x <C_\beta
\end{equation}
where
$$C_\alpha=\frac{8\sin(x/2)-\sin(x)}{\alpha} \quad {\rm and}\quad C_\beta=\frac{8\sin(x/2)-\sin(x)}{\beta},$$
with the best possible constants $\alpha=3$ and $\beta=(8\sqrt{2}-2)/\pi \approx 2.96465$.
\end{theorem}

By using Mathematica Software\textsuperscript{\tiny{\textregistered}} \cite{ru}, one can see that
Theorem \ref{thm1} refines the
inequalities in \eqref{zhuineq} as follows:
$$Z_l(x)<C_\alpha(x),\quad {\rm for}\quad x\in(0,1.28966),$$
$$Z_u(x)<C_\beta(x),\quad {\rm for}\quad x\in(0,0.980316).$$
It is worth to mention that the first inequality in Theorem \ref{thm1} 
was discovered heuristically by Huygens \cite{campan}, here we have given a
proof.

In 1970, Carlson \cite{carl} established the following inequalities,
\begin{equation}\label{carlineq}
\frac{6(1-x)^{1/2}}{2\sqrt{2}+(1+x)^{1/2}}<\arccos(x)<\frac{4^{1/3}(1-x)^{1/2}}{(1+x)^{1/6}},
\end{equation}
$0<x<1$. These inequalities are known as
Carlson's inequalities in the literature. Thereafter, several authors studied these
inequalities, and gave some generalization and partial refinement,
e.g., see \cite{chen,guoqi,{guoqi2},z3}. It is interesting to
observe that the Adamovi\'c-Mitrinovi\'c and Cusa-Huygens inequality
\eqref{cusa} implies the second and the first inequality of
\eqref{carlineq}, respectively, with the transformation
$x=\arccos(t),\,0<t<\pi/4$.

For $0<x<1$, Guo and Qi \cite{guoqi,guoqi2} gave the following inequalities,
\begin{equation}\label{guoineq1}
\frac{\pi}{2}\frac{(1-x)^{1/2}}{(1+x)^{1/6}}<\arccos(x)
<\frac{(1/2+\sqrt{2})(1-x)^{1/2}}{2\sqrt{2}+(1+x)^{1/2}},
\end{equation}

\begin{equation}\label{guoineq2}
\frac{4^{1/\pi}(1-x)^{1/2}}{(1+x)^{(4-\pi)/(2\pi)}}<\arccos(x)<
\frac{\pi(1-x)^{1/2}}{2(1+x)^{(4-\pi)/(2\pi)}}.
\end{equation}
They concluded that these inequalities doesn't refine \eqref{carlineq} in the whole interval $(0,1)$ of $x$.

Chen et al. \cite{chen} established the lower bound for $\arccos$ as follows,
\begin{equation}\label{chenineq}
\frac{\pi}{2}\frac{(1-x)^{(\pi+2)/\pi^2}}{(1+x)^{(\pi-2)/\pi^2}}<\arccos(x),\quad 0<x<1.
\end{equation}
The inequality \eqref{chenineq} refines the first inequality of \eqref{carlineq} for $x\in(0,0.345693)$.

In \cite{z3}, Zhu proved that for $p\geq 1$ and $x\in(0,1)$
\begin{equation}\label{z3ineq}
\frac{2\cdot 3^{1/p}\sqrt{1-x}}{((2\sqrt{2})^p+(\sqrt{1+x})^p)^{1/p}}<\arccos(x)<
\frac{2\pi\sqrt{1-x}}{((2\sqrt{2})^p+(\pi^p-2^p)(\sqrt{1+x})^p)^{1/p}},
\end{equation}
inequalities reverse for $p\in[0,4/5]$.

We give the following theorem, which refines the Carlson's inequality, see Figure 1.

\begin{theorem}\label{thm3} For $x\in(0,1)$
\begin{equation}\label{thm3inequ}
\frac{1}{3}\left(8\sqrt{2-\sqrt{2}\sqrt{1+x}}\right) < \arccos(x) <
\frac{2^{11/6}\sqrt{1-x}}{(2+\sqrt{2}\sqrt{1+x})^{2/3}}.
\end{equation}
\end{theorem}

We see that Theorem \ref{thm3} refines the inequalities in \eqref{z3ineq} by using the Mathematica Software\textsuperscript{\tiny{\textregistered}} \cite{ru}.

In 1967, Shafer \cite{sha} proposed the following elementary inequality
\begin{equation}\label{shaineq}
\frac{3x}{1+2\sqrt{1+x^2}}<\arctan(x),\quad x>0.
\end{equation}
This inequality was proved by Grinstein, Marsh and Konhauser by different ways in \cite{sha2}.

In 2009, Qi et al. \cite{qizha} refined the inequality \eqref{shaineq} as follows,

\begin{equation}\label{qiineq}
\frac{(1+a)x}{a+\sqrt{1+x^2}}<\arctan(x)<\frac{(\pi/2)x}{4+\sqrt{1+x^2}},\,x>0,\,-1<a<1/2,
\end{equation}
$$\frac{4a(1+a^2)x}{a+\sqrt{1+x^2}}<\arctan(x)<\frac{\max\{\pi/2,1+a\}x}{a+\sqrt{1+x^2}},\,x>0,\,1/2<a<2/\pi.$$

Recently, Alirezaei \cite{ali} has sharpened the Shafer's inequality \eqref{shaineq} by giving the following bounds for
$\arctan$,
\begin{equation}\label{aliinequ}
\frac{x}{4/\pi^2+\sqrt{(1-4/\pi^2)^2+4x^2/\pi^2}} <\arctan(x)
\end{equation}
$$<\frac{x}{1-6/\pi^2+\sqrt{(6/\pi^2)^2+4x^2/\pi^2}},$$
for $x\in\mathbb{R}$.
Graphically, it is shown that the maximum relative errors of
the obtained bounds are approximately smaller than $0.27\%$ and
$0.23\%$ for the lower and upper bound, respectively.

Our next result refines the bounds given in \eqref{aliinequ}, which is illustrated in Figure 2.

\begin{theorem}\label{thm4} For $x\in(0,1)$, we have
\begin{equation}\label{thm4inequ}
\frac{1}{3}\left(4\sqrt{2}\sqrt{1-\frac{1}{\sqrt{1+x^2}}}-\frac{x}{\sqrt{1+x^2}}\right)<\arctan(x)
\end{equation}
$$<\frac{2^{2/3} x}{\sqrt{1+x^2}
   \left(1+1/\sqrt{1+x^2}\right)^{2/3}}.$$
\end{theorem}



\section{Preliminaries}

For easy reference, we recall the the following Monotone l'H\^opital
rule due to Anderson et al. \cite[Theorem 2]{avv1}, which has been
extremely used in the literature.
\begin{lemma}\label{lem0}
For $-\infty<a<b<\infty$,
let $f,g:[a,b]\to \mathbb{R}$
be continuous on $[a,b]$, and be differentiable on
$(a,b)$. Let $g^{'}(x)\neq 0$
on $(a,b)$. If $f^{'}(x)/g^{'}(x)$ is increasing
(decreasing) on $(a,b)$, then so are
$$\frac{f(x)-f(a)}{g(x)-g(a)}\quad and \quad \frac{f(x)-f(b)}{g(x)-g(b)}.$$
If $f^{'}(x)/g^{'}(x)$ is strictly monotone,
then the monotonicity in the conclusion
is also strict.
\end{lemma}

\begin{lemma}\label{lemma1702a} The function
$$f(x)=4x\sin(x)+(4-x^2)\cos(x)-x^2$$
is strictly decreasing from $(0,\pi/2)$ onto $(a,4),\,a=\pi(8-\pi)/4\approx 3.81578$. In particular,
$$\frac{\pi(8-\pi)/4+x^2-(4-x^2)\cos(x)}{4x^2}<\frac{\sin(x)}{x}<\frac{4+x^2-(4-x^2)\cos(x)}{4x^2}$$
for $x\in(0,\pi/2)$.
\end{lemma}

\begin{proof} By differentiating and using the indentities
$\sin(x)=2\sin(x/2)\cos(x/2)$ and $1-\cos(x)=2\sin(x/2)^2$ we get
\begin{eqnarray*}
f'(x)&=& x(2\cos(x)+x\sin(x)-2)\\
&=&2\sin(x/2)(x\cos(x/2)-2\sin(x/2))<0.
\end{eqnarray*}
Hence $f$ is strictly decreasing in $x\in(0,\pi/2)$, and the limiting values can be obtained easily.
\end{proof}

\begin{lemma}\label{lemma1702b} The following function
$$f(x)=\frac{\sin(x)-x\cos(x)}{2\sin(x/2)-x\cos(x/2)}$$
is strictly decreasing from $(0,\pi/2)$ onto $(b,4),\,b=2\sqrt{2}/(4-\pi)\approx 3.81578$. In particular,
$$\frac{2\sqrt{2}}{4-\pi}\left(2\sin\left(\frac{x}{2}\right)-x\cos\left(\frac{x}{2}\right)\right)<\sin(x)-x\cos(x)<
4\left(2\sin\left(\frac{x}{2}\right)-x\cos\left(\frac{x}{2}\right)\right),$$
for $x\in(0,\pi/2)$.
\end{lemma}

\begin{proof} We get
\begin{eqnarray*}
f'(x)&=&\frac{x\sin(x)}{2\sin(x/2)-x\cos(x/2)}-\frac{x\sin(x/2)(\sin(x)-x\cos(x))}{2(2\sin(x/2)-x\cos(x/2))^2}\\
&=&\frac{x\sin(x/2)(x(2+\cos(x))-3\sin(x))}{4-(4x\sin(x)+(4-x^2)\cos(x)-x^2)},
\end{eqnarray*}
which is negative by the second inequality of \eqref{cusa} and Lemma \ref{lemma1702a}.
This implies that $f$ is strictly
decreasing $x\in(0,\pi/2)$, and by applying l'H\^opital rule we get the limiting values.
\end{proof}

%
%


\begin{lemma}\label{1702c} The following function
$$g(x)=\frac{8\sin(x/2)-\sin(x)}{x}$$
is strictly decreasing from $(0,\pi/2)$ onto $(\alpha,3),\,\beta=(8\sqrt{2}-2)/\pi \approx 2.96465$. Also, the function
$$f(z)=\frac{8\sin(z)}{6z+\sin(2z)}$$
is strictly decreasing from $(0,\pi/4)$ onto $(1,\gamma),\,\gamma=8\sqrt{2}/(2+3\pi)\approx 0.99028$.
\end{lemma}

\begin{proof} We get
\begin{eqnarray*}
g'(x)&=&\frac{4\cos(x/2)-\cos(x)}{x}-\frac{8\sin(x/2)-\sin(x)}{x^2}\\
&=&\frac{\sin(x)-x\cos(x)-4(2\sin(x/2)-x\cos(x/2))}{x^2},
\end{eqnarray*}
which is negative by Lemma \ref{lemma1702b}. Thus, $g$ is strictly
decreasing in $x\in(0,\pi/2)$ , and the limiting values follow from the l'H\^opital
rule.

Next, let
$f=f_1(z)/f_2(z),\,z\in(0,\pi/4)$, where $f_1(z)=8\sin(z)$ and $f_2(z)=6z+\sin(2z)$.
We get
$$\frac{f'_1(z)}{f'_2(z)}=\frac{4\cos(z)}{1+\cos(z)^2}=f_3(z).$$
One has,
$$f'_3(z)=-\frac{\sin(z)^3}{(3+\cos(2z))^2}<0.$$
Clearly, $f_1(0)=f_2(0)=0$, hence by Lemma \ref{lem0} $f$ is strictly decreasing, and we get
$$\lim_{z\to \pi/4}f(z)=8\sqrt{2}/(2+3\pi)\approx 0.99028 <f(z)< \lim_{z\to 0}f(z)=1.$$
This implies the proof if we let $z=x/2$.
\end{proof}


\section{Proof of Theorems}

\vspace{.3cm}
\noindent{\bf Proof of Theorem \ref{thm1}.} \rm The proof follows easily
from Lemma \ref{1702c}. $\hfill\square$

\begin{corollary} For $x\in(0,\pi/2)$, we have
$$
\frac{8\sin(x/2)-\sin(x)}{\beta}
<\frac{8\sin(x/2)-\beta \sin(x)}{\gamma},
$$
where $\beta$ and $\gamma$ are as in Lemma \ref{1702c}.
\end{corollary}

\begin{proof}
For $x\in(0,\pi/2)$, let
$f(x)=\sin(x/2)/\sin(x)$. One has
$$f'(x)=\frac{\sin(x/2)^3}{\sin(x)^2}>0.$$
Hence, $f$ is strictly increasing, and
$$\frac{1}{2}=\lim_{x\to 0}f(x)<f(x)<\lim_{x\to \pi/2}f(x)=\frac{1}{\sqrt{2}}.$$
We observe that
$$\frac{\sin(x/2)}{\sin(x)}<\frac{1}{\sqrt{2}}=-\frac{2-8\sqrt{2}+3\pi}{16-2\sqrt{2}-3\sqrt{2}\pi},$$
which is equivalent to
$$\frac{(16-2\sqrt{2}-3\sqrt{2}\pi)\sin(x/2)+(2-8\sqrt{2}+3\pi)\sin(x)}{24\sqrt{2}-6}>0.$$
This is equivalent to the desired inequality.
\end{proof}

\vspace{.3cm}
\noindent{\bf Proof of Theorem \ref{thm3}.} \rm
Let $x=\cos(2t)$ for $0<t<\pi/4$. Then $\arccos(x)/2=t$, and clearly $0<x<1$.
From \eqref{refine} and \eqref{thm1inequ} we have
\begin{equation}\label{1802a}
\frac{8\sin(t/2)-\sin(t)}{3}<t<\frac{2^{2/3}\sin(t)}{(1+\cos(t))^{2/3}},
\end{equation}
for $t\in(0,\pi/2)$. Replacing $\cos(t),\,\sin(t)$ and $t$ by
$\sqrt{(1+x)/2},\,\sqrt{(1-x)/2}$ and $\arccos(t)/2$, respectively,
in \eqref{1802a}, we get
$$\frac{8((1-\sqrt{(1+x)/2})/2)^{1/2}-\sqrt{(1-x)/2}}{3}<\frac{\arccos(x)}{2}<
\frac{2^{2/3}\sqrt{(1-x)/2}}{(1+\sqrt{(1+x)/2})^{2/3}}.$$ After
implication we get the desired inequality.$\hfill\square$

\vspace{.3cm}
\noindent{\bf Proof of Theorem \ref{thm4}.} \rm Next, let
$x=\tan(t), \ t\in(0,\pi/2)$ and $x\in(0,1)$. Then
$t=\arctan(x)$, and by using the identity $1+\tan(t)^2=\sec(t)^2$ we get
$$\sin(t)=\frac{x}{\sqrt{1+x^2}}=m,\quad {\rm and}\quad \sin\left(\frac{t}{2}\right)=
\left(\frac{\sqrt{1+x^2}-1}{2\sqrt{1+x^2}}\right)=n.$$
We get the desired inequalities if we replace, $t,\,\sin(t),\,\sin(t/2)$ by $\arctan(x),\,m,\,n$, respectively,
in \eqref{1802a}. $\hfill\square$

\vspace{.3cm}
For the comparison of the bounds of $\arccos(x)$ and $\arctan(x)$ given in \eqref{carlineq} and \eqref{thm3inequ}
with the corresponding bounds appear in Theorem \ref{thm3} and \ref{thm4}, we use the the graphical method,
see Figure 1 and 2.

\begin{figure}\label{arccosfi}
\includegraphics[width=12cm]{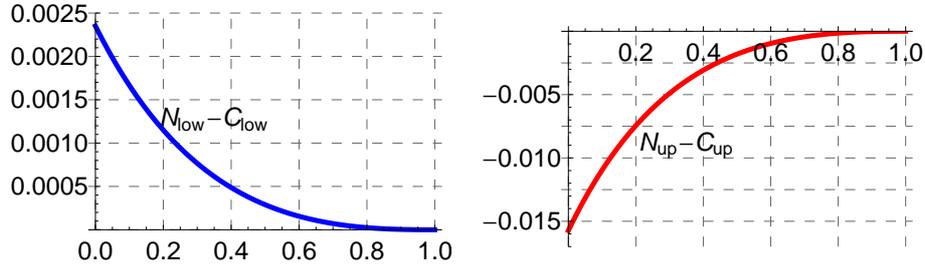}
\caption{We denote the left-hand sides of \eqref{carlineq} and \eqref{thm3inequ} by $C_{low}$ and $N_{low}$, respectively, while the right-hand sides by $C_{up}$ and $N_{up}$, respectively. It is clear that
\eqref{thm3inequ} refines the Carlson's inequality \eqref{carlineq}.}
\end{figure}

\begin{figure}
\includegraphics[width=12cm]{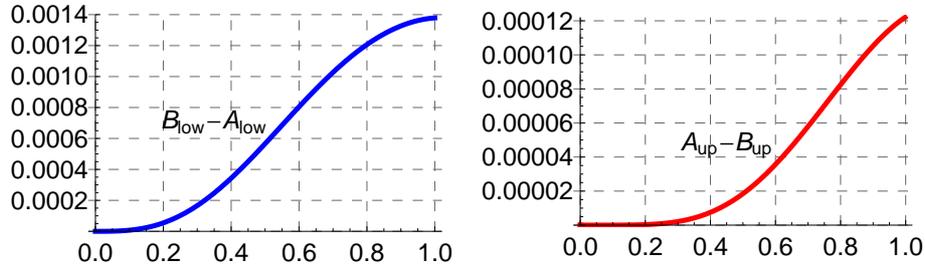}
\caption{We denote the lower and upper bound of \eqref{thm4inequ} by
$B_{low}$ and $B_{up}$, respectively,  while the corresponding
bounds of \eqref{aliinequ} are denoted by $A_{low}$ and $A_{up}$.
The differences $B_{low}-A_{low},\ A_{up}-B_{up}$ are positive, this
implies that the inequalities in \eqref{thm4inequ} are better than
the corresponding inequalities of \eqref{aliinequ}.}
\end{figure}

\vspace{.3cm}
The proof of the following corollary is the analogue of the proofs of the
previous theorems, and can be obtained by utilizing \eqref{zhuineq}.

\begin{corollary}\label{coro1} For $x\in(0,1)$, we have
\begin{enumerate}
\item $\displaystyle\frac{\pi\sqrt{1-x}}{\sqrt{2}+(2/\pi)\sqrt{1+x}}
<\arccos(x)<\frac{\pi\sqrt{1-x}}{\sqrt{2}+(\pi/2-1)\sqrt{1+x}},$\\
\item $\displaystyle\frac{\pi/2x}{2/\pi+\sqrt{1+x^2}}<\arctan(x)<\frac{\pi x}{(\pi-2)+\sqrt{1+x^2}}.$
\end{enumerate}
\end{corollary}



%
%


\begin{abstract}
In this paper, authors refine the classical inequalities of the trigonometric functions, such as Jordan's inequality, Cusa-Huygens inequality and Kober's inequality.
\end{abstract}

\maketitle

%
%


\section{Introduction}
The study of the classical inequalities of the trigonometric functions such as
Adamovi\'c-Mitrinovi\'c inequality, Cusa-Huygens inequality, Jordan inequality, Redheffer inequality,
Becker-Stark inequality, Wilker inequality,
Huygens inequality, and Kober inequality has got big attention of the numerous authors. Since last ten years, the huge number of papers on the refinement and the generalization of these inequalities have appeared, e.g. see
\cite{avv1,avv2,kvv,lili,neua,neu2a,neusana,joz2,joz3,sanben,yang} and the references therein. Motivated by these studies, 
in this paper we refine Jordan's, Kober's and Cusa-Huygens inequality, and our results refine the existing results in the literature.



%


%

The well-know Jordan's inequality \cite{mit} states,
\begin{equation}\label{jordan}
\frac{\pi}{2}\leq \frac{\sin(x)}{x},\quad 0<x\leq \frac{\pi}{2},
\end{equation}
equality with $x=\pi/2$.

In 2003, Debnath and Zhao \cite{deb} refined the inequality \eqref{jordan} as below,
\begin{equation}\label{debinequ1}
d_1(x)=\frac{2}{\pi }+\frac{1}{12\pi}(\pi ^2-4 x^2)\leq\frac{\sin(x)}{x},
  \end{equation}
    \begin{equation}\label{debinequ2}
d_2(x)=\frac{2}{\pi }+\frac{1}{\pi^3}(\pi ^2-4 x^2)\leq\frac{\sin(x)}{x},
  \end{equation}
    for $x\in(0,\pi/2)$, equality in both inequalities with $x=\pi/2$. 
Thereafter, another proof of the inequality \eqref{debinequ2} was given by Zhu in \cite{z2}.

In 2006,  \"Ozban \cite{ozban} proved the following inequality,
    \begin{equation}\label{ozbaninequ1}
 o(x)=\frac{2}{\pi }+\frac{1}{\pi^3}(\pi ^2-4 x^2)+\frac{4(\pi-3)}{\pi^3}\left(x-\frac{\pi}{2}\right)^2\leq
\frac{\sin(x)}{x},
  \end{equation}
     for $x\in(0,\pi/2)$, equality with $x=\pi/2$. 
%

In the same year, the following refinement of \eqref{jordan} was proved by Jiang and Yun \cite{jiang},
\begin{equation}\label{doninequ}
j(x)=\frac{2}{\pi }+\frac{\pi ^4-16 x^4}{2\pi^5}<\frac{\sin(x)}{x},\quad  0<x<\frac{\pi}{2}.
\end{equation}
for $x\in(0,\pi/2)$, equality with $x=\pi/2$. 

In \cite{zhang3003}, Zhang et al. gave the following inequality,
\begin{equation}\label{ineq3003a}
zw(x)=\frac{3}{\pi}-\frac{4}{\pi^3}x^2<\frac{\sin(x)}{x},\quad 0<x<\frac{\pi}{2}.
\end{equation}

It is easy to see that $d_1(x)<d_2(x)$, $d_2(x)=zw(x)$ and $j(x)<d_2(x)<o(x)$ for $x\in(0,\pi/2)$.


Our first main result refines the inequality \eqref{ozbaninequ1} as follows:

\begin{theorem}\label{thm2} For $x\in(0,\pi/2)$, we have
$$o(z)\leq 1+\frac{16(\pi-3)}{\pi^4}x^3-\frac{4(3\pi-8)}{\pi^3}x^2\leq\frac{\sin(x)}{x},$$
equality in both inequalities holds with $x=\pi/2$.
\end{theorem}

In literature, the following inequalities
\begin{equation}\label{laztri}
(\cos x)^{1/3}<\frac{\sin x}{x}<\frac{\cos x+2}{3},\quad 0<|x|<\frac{\pi}{2},
\end{equation}
are known as Adamovi\'c-Mitrinovi\'c inequality \cite[p.238]{mit} and Cusa-Huygens \cite{sanben} inequality, respectively.
For the refinement of \eqref{laztri}, e.g. see \cite{kvv,neua,neusana,joz3,sanben,yang} and the bibliography of these papers. Most of the refinements of \eqref{laztri} are involving very complicated upper and lower bound of $\sin(x)/x$. In the following theorem we refine  \eqref{laztri} by giving the upper and lower bound of $\sin(x)/x$ in terms of simple functions, and these functions are also independent of the exponent.

\begin{theorem}\label{thm1} For $x\in(0,\pi/2)$, we have
$$\frac{1+\cos(x)}{2-\alpha x^2}<\frac{\sin(x)}{x}<\frac{1+\cos(x)}{2-\beta x^2},$$ 
with the best possible constants $\alpha=1/6\approx 0.166667$ and $\beta=2(4-\pi)/\pi^2\approx 0.17396$.
\end{theorem}

In 1944, Kober \cite[3.4.9]{mit} established the following inequalities,
\begin{equation}\label{koberineq}
1-2\frac{x}{\pi}<\cos(x)<1-2\frac{x}{\pi},
\end{equation}
the first inequality is valid for $x\in(0,\pi/2)$, while the second one holds for  $x\in(\pi/2,\pi)$. In literature, these
inequalities are known as Kober's inequalities.

By studying the function $x\mapsto (1-\cos(x))/x,\,x\in(0,\pi/2)$,  S\'andor \cite{sandor17b} refined the Kober's inequality as follows:
\begin{equation}\label{ineq3103a}
\cos(x)<1-\frac{2}{\pi}x-\frac{2(\pi-2)}{\pi^2}\left(x-\frac{\pi}{2}\right), \quad 0<x<\frac{\pi}{2},
\end{equation}
\begin{equation}\label{ineq3103b}
1-\frac{x^2}{2}<\cos(x)<1-\frac{4x^2}{\pi^2}, \quad 0<x<\frac{\pi}{2}.
\end{equation}

In \cite{zhang3003}, the following refinement appeared,
\begin{equation}\label{ineq3003}
1-\frac{4-\pi}{\pi}x-\frac{2(\pi-2)}{\pi^2}x^2<\cos(x)<1-\frac{4}{\pi^2}x^2,\quad 0<x<\frac{\pi}{2}.
\end{equation}

Applying Taylor series expansion one can have
\begin{equation}\label{ineq3103e}
1-\frac{x^2}{2}<\cos(x)<1-\frac{x^2}{2}+\frac{x^4}{24},\quad 0<x<\frac{\pi}{2}.
\end{equation}

On the basis of Mathematica Software\textsuperscript{\tiny{\textregistered}} \cite{ru}, we conclude that our next result refines the above Kober's inequalities as follows:
\begin{theorem}\label{prop1702} For $x\in(0,\pi/2)$, we have
$$1-\frac{x^2/2}{1+x^2/12}<\cos(x)<1-\frac{24x^2/(5\pi^2)}{1+4x^2/(5\pi^2)}.$$
\end{theorem}

\section{Preliminaries and proof of main result}

The following lemma will be used in the proof of theorems.

\begin{lemma}\label{lema1} (1) The function 
$$f_1(x)=\cos(x)^{2/3}+\frac{x^2}{3}$$
is strictly decreasing from $(0,\pi/2)$ onto $(\pi^2/12,1)$. In particular, for $x\in(0,\pi/2)$
$$\left(\frac{\pi^2-4x^2}{12}\right)^{3/2}<\cos(x)<\left(1-\frac{x^2}{3}\right)^{3/2},$$
equivalently,
$$\frac{x^2\cos(x)^{1/3}+3\cos(x)}{3}<\cos(x)^{1/3}<\frac{4}{\pi^2}\left(\frac{x^2\cos(x)^{1/3}+3\cos(x)}{3}\right).$$
(2) The function
$$f_2(x)=1+\cos(x)^{4/3}-2\,\frac{\sin(x)^2}{x^2}$$
is strictly increasing from $(0,\pi/2)$ onto $(0,1-8/\pi^2)$. In particular, for $x\in(0,\pi/2)$
$$\frac{\cos(x)^{4/3}+8/\pi^2}{2}<\frac{\sin(x)^2}{x^2}<\frac{1+\cos(x)^{4/3}}{2},$$
and 
$$\frac{4(1-\cos(x))}{x^2}-1<\cos(x)^{4/3}<\frac{4(1-\cos(x))}{x^2}+\frac{1}{2^{2/3}}-\frac{16}{\pi^2}.$$
\end{lemma}

\begin{proof} By Adamovi\'c-Mitrinovi\'c inequality, we get
$$f_1'(x)=\frac{2x}{3}-\frac{2\sin(x)}{3\cos(x)^{1/3}}=
\frac{2x}{3}\left(1-\frac{\sin(x)/x}{\cos(x)^{1/3}}\right)<0.$$
Since, $f_1$ is strictly decreasing in $x\in(0,\pi/2)$, and the limiting values are clear.

Next for the proof of part (2), we get
$$f_2'(x)=-\frac{4\sin(x)}{x^3}\left(\frac{x^3\cos(x)^{1/3}}{3}+x\cos(x)-\sin(x)\right)>
-\frac{4\sin(x)}{x^3}\left(\cos(x)^{1/3}-\sin(x)/x\right)>0,$$
by part (1).
Thus, $f_2$ is strictly increasing in $x\in(0,\pi/2)$, and the limiting values follow easily.
Last inequalities follow if we consider $f_2(x/2),x\in(x,\pi/2)$.
\end{proof}


\vspace{.3cm}
\noindent{\bf Proof of Theorem \ref{thm1}.} For $x\in(0,\pi/2)$, let 
$$g_1(x)=\frac{(1+\cos(x))}{x\sin(x)}-\frac{2}{x^2}.$$
By Lemma \ref{lema1}(2) we get
\begin{eqnarray*}
g_1'(x)&=&\frac{4}{x^3}-\frac{1}{x}-\frac{1+\cos(x)}{x^2\sin(x)}-\frac{(1+\cos(x))\cos(x)}{x\sin(x)}\\
&=& \frac{4(1-\cos(x))/x^2-1-\sin(x)/x}{x(1+\cos(x))}< \frac{\cos(x)^{4/3}-\sin(x)/x}{x(1+\cos(x))}<0.
\end{eqnarray*}
Since $g_1$ is strictly decreasing in $x\in(0,\pi/2)$, and by applying l'H\^opital rule, we get the liming values.
This completes the proof.$\hfill\square$

\begin{lemma}\label{lema2} The function
$$f_3(x)=\frac{x-\sin(x)}{x^3}$$
is strictly decreasing and concave from $(0,\pi/2)$ onto $(1/\pi^2,1/6)$. In particular, for $x\in(0,\pi/2)$
$$1-\frac{x^2}{6}<\frac{\sin(x)}{x}<1-\frac{x^2}{\pi^2}.$$
\end{lemma}

\begin{proof} One has,
$$f_3'(x)=\frac{1-\cos(x)}{x^3}-3\frac{x-\sin(x)}{x^4}=\frac{3\sin(x)/x-2+\sin(x)}{x^4},$$
which is negative by Cusa-Huygens inequality, hence $f_3$ is strictly decreasing in $x\in(0,\pi/2)$. Again
\begin{eqnarray*}
f_3''(x)&=&\frac{2\cos(x)+x\sin(x)-2}{x^5}-\frac{4}{x^5}(3\sin(x)-x(2+\cos(x)))\\
&=& \frac{6x(1+\cos(x))-(12-x^2)\sin(x)}{x^5},
\end{eqnarray*}
which is negative by Theorem $\ref{thm1}$, this implies the concavity of the function $f_3$.
\end{proof}

\noindent{\bf Proof of Theorem \ref{thm2}.}
Using the concavity of the function $f_3(x),\,x\in(0,\pi/2)$, the tangent line on point $(\pi/2, f_3(\pi/2))$ is above the graph of $f_3(x)$ on $(0, \pi/2)$. The equation of the tangent line is 
$$y=\frac{4(\pi-2)}{\pi^3}+\frac{16(3-\pi)}{\pi^4}(x-\pi/2).$$
After some computations, we get the desired inequality.
The first inequality is equivalent to 
$$-\frac{4 (\pi -3) (\pi -2 x)^2
   x}{\pi ^4}<0,$$
	which is obvious. This completes the proof.$\hfill\square$

%
%
\vspace{.3cm}
\noindent{\bf Proof of Theorem \ref{prop1702}.} For $x\in(0,\pi/2)$, let
$$f(x)=\frac{x^2(5+\cos(x))}{1-\cos(x)}.$$
Differentiation gives
$$f'(x)=\frac{2x(5-g(x))}{(\cos(x)-1)^2},$$
where
$$g(x)=\cos(x)(4+\cos(x))+3x\sin(x).$$
We get
\begin{eqnarray*}
g'(x)&=&3x\cos(x)-(1+2\cos(x))\sin(x)\\
&=&x\cos(x)\left(3-\left(2\frac{\sin(x)}{x}+\frac{\tan(x)}{x}\right)\right),
\end{eqnarray*}
which is negative by so-called Huygens inequality \cite{neusana} 
$$2\frac{\sin(x)}{x}+\frac{\tan(x)}{x}>3,\quad 0<x<\frac{\pi}{2}.$$
Thus, $g$ is decreasing and $\lim_{x\to 0}g(x)=5,$
and as a result $f'>0$. This implies that $f$ is strictly increasing. 
Applying l'H\^opital rule we get
$$12=\lim_{x\to 0}f(x)<f(x)<\lim_{x\to 0}f(x)=\frac{5\pi^2}{4}\approx 12.33701,$$
which is equivalent to 
$$\frac{6}{1+x^2/12}-5<\cos(x)< \frac{6}{1+4x^2/(5\pi^2)}-5.$$
This implies the desired inequalities. $\hfill\square$




%


\begin{abstract}
We offer two new proofs of famous Redheffer's
inequality, as well establish two converse inequalities for it. 
Also a hyperbolic analogue is pointed out.
\end{abstract}

\maketitle
\bigskip
{\bf 2010 Mathematics Subject Classification}: 26D05, 26D07, 26D99.

{\bf Keywords}: Trigonometric and hyperbolic functions, Redheffer's inequality, Cusa-Huygens inequality.


\section{Introduction}
In 1969, Redheffer \cite{re} proposed the following inequality
\begin{equation}\label{sinRed}
\frac{\pi^2-x^2}{\pi^2+x^2}\leq \frac{\sin x}{x},\quad x\in\mathbb{R},
\end{equation}
which was proved by Williams \cite{wi} in the same year.
 In literature, this inequality is known as
\emph{Redheffer's inequality}. By using the infinite product and induction
method, Williams verified this inequality also in 1969, \cite{wi1}. Motivated by his work, many developments such as
generalizations, refinements and applications took place, e.g., see \cite{kua}.

Thereafter some Redheffer's-type inequalities for other trigonometric, hyperbolic and
Bessels function were established, e.g., see \cite{bar,czq,z1,sz} and the references therein.


Recently a new proof of \eqref{sinRed} has appeared in \cite{lz}, where authors
are using the Lagrange mean value theorem, combined with  induction , which is very complicated for the reader.

The inequality \eqref{sinRed} is valid for all $x\in\mathbb{R}$.
It is immediate that we may assume $x>0$ and $x\in(0,\pi)$ as for $x>\pi$,
we may let $x=\pi+t$ for $t>0$,
then inequality \eqref{sinRed} becomes
$$\frac{\sin t}{t} < \frac{2\pi^2+3\pi t+t^2}{2\pi^2+2\pi t+ t^2}.$$
This is obvious, as right side is greater than one, and left side less than one.
Thus, we may consider $x \in (0,\pi)$.

So far, all the authors have given the proof of \eqref{sinRed} by
using the induction method. From our proof, it is obvious that the
induction method is not needed. In this paper we give new interesting proofs of
\eqref{sinRed}, which are based on the elementary calculus.  The
authors think that this proof could be one from the ``Book" (See
[2]).

\section{New proof of inequality and its converse}

\begin{lemma}\label{lemma1}
For all $x\in(0,\pi)$ one has
$$x+ \sin( x)  > x^2\frac{\cos (x/2)}{\sin (x/2)}.$$
\end{lemma}

\begin{proof} Let $h_1(x)= x+\sin (x)-x^2\cos (x/2)/\sin (x/2).$
Then
$$h_1'(x)= 1+ \cos (x) -\frac{2x\sin (x)-x^2}{2\sin (x/2)^2}=
\frac{\sin (x)^2 +x^2-2x\sin(x)}{2\sin(x/2)^2} >0, $$
 where we have used $ 1+\cos( x)=2\cos (x/2)^2$  and $2\sin(x/2)\cos(x/2)=\sin(x)$. As
$h_1(x)  >h_1(0)=0$, the result follows.
\end{proof}

\begin{theorem}\label{thm1} The following inequalities
\begin{equation}\label{ineq}
\frac{\pi^2-x^2}{\pi^2+x^2}\leq \frac{\sin(x)}{x}< \frac{12-x^2}{12+x^2}
\end{equation}
hold for  $x \in(0,\pi]$. 
\end{theorem}

\begin{proof} Let
$$f_1(x)=\frac{x^2(x+\sin( x))}{x-\sin( x)},$$
for $x \in(0,\pi]$.
After some elementary computations, one has
$$\frac{(x-\sin (x))^2}{2x}f'(x)=  x^2(1+\cos(x))-\sin(x)(x+\sin (x))= g_1(x).$$
As $1+ \cos(x)=2\cos(x/2)^2$, and $\sin (x)=2\sin (x/2)\cos(x/2)$,
we get
\begin{eqnarray*}
g_1'(x)&=& 2\cos \left(\frac{x}{2}\right)\left[x^2\cos \left(\frac{x}{2}\right)-
\sin \left(\frac{x}{2}\right)(x+\sin (x))\right]\\
&=&2\sin \left(\frac{x}{2}\right)\cos\left(\frac{x}{2}\right)
\left[x^2\frac{\cos(x/2)}{\sin(x/2)}-(x+\sin(x))\right]<0
\end{eqnarray*}
by Lemma \ref{lemma1}. Thus $f_1'(x)<0$, and $f_!$ is strictly decreasing in $x\in(0,\pi]$.
We get $f_1(x)>\lim_{x\to \pi}f_1(x)=\pi^2$, which is equivalent to 
$$\frac{x+\sin(x)}{x-\sin(x)} > \frac{\pi^2}{x^2},$$
thus the first inequality in \eqref{ineq} follows. 
The second inequality in \eqref{ineq} follows similarly from $f_1(x)<\lim_{x\to 0 }f_1(x)=12$.
\end{proof}

%

The right side of \eqref{ineq} improves the known inequality:
$$\frac{\sin (x)}{x}< 1-\frac{x^2}{\pi^2},\quad 0<x<\frac{\pi}{2}.$$
Indeed, this follows by $x^2 +12 <2\pi^2$, which is true, as
$\pi^2/4+12<2\pi^2$  becomes $48<7\pi^2$.


\begin{proposition}
The second inequality in \eqref{ineq} refines the relation
$$\frac{\sin(x)}{x}<\frac{\cos(x)+2}{3},\quad 0<x<\frac{\pi}{2},$$
so-called Cusa-Huygens inequality \cite{neusan,sanben}.
\end{proposition}

\begin{proof}
It is equivalent to prove that,
$$(12+x^2)\cos(x)+5x^2-12=s(x)>0.$$
One has successively:
$$s'(x)= 2x(5+\cos (x))-(12+x^2)\sin(x) -12,$$
$$s''(x)= -4x\sin (x)-(10 +x^2)\cos (x) +10,$$
$$s'''(x)= x^2\sin(x) +6(\sin(x)-x\cos(x))>0,$$
 as $\sin(x)>x\cos(x)$ (i.e. $\tan(x) >x$).
Thus we get  $$s''(x)>s''(0)=0,\, s'(x)>s'(0)=0,$$
and finally $s(x)>s(0)=0$.
\end{proof}



\section{An other proof of $\eqref{sinRed}$}

\begin{lemma}\label{newlema}
Let $g(x)=f(x)/\sin (x)$ for $x \in (0, \pi)$. Then $\sin(x)^2g'(x)= h(x)$, and  the
sign of $h'(x)$ depends on the sign of $F(x)= f(x)+ f''(x)$.
\end{lemma}

\begin{proof} One has  
$$\sin (x)^2 g'(x)= f'(x)\sin(x)-f(x)\cos (x)=h(x),$$ and 
$$h'(x)= (f(x)+f''(x))\sin(x)=F(x)\sin(x).$$
As $\sin (x)>0$ for all $x \in (0,\pi)$, the result follows.
\end{proof}

\begin{theorem}\label{thm16} For $x\in(0,\pi)$, we have
$$\frac{\pi^2-x^2}{\pi^2+x^2}<\frac{\sin(x)}{x}< c_1\frac{\pi^2-x^2}{\pi^2+x^2},$$
where $c_1=1.07514$.
\end{theorem}

\begin{proof} Let $g(x)=f(x)/\sin (x)$ for $x \in (0, \pi)$, where
$$f(x)=\frac{x(\pi^2-x^2)}{\pi^2+x^2}.$$
One has 
$$ f'(x)=  \frac{\pi^4-4\pi^2x^2-x^4}{(\pi^2+x^2)^2},$$  so we get
$$F(x)= f(x)+ f''(x)= \frac{ x p(x)}{(\pi^2+x^2)^3},$$
i.e, by Lemma \ref{newlema} the sign of  $h'(x)$ depends on the sign of $p(x)$.
Here 
\begin{eqnarray*}
h(x)&=& f'(x)\sin(x)-f(x)\cos(x)\\
&=&\frac{x \left(\pi ^4-x^4\right) \cos (x)+
\left(-x^4-4 \pi ^2 x^2+\pi ^4\right) \sin (x)}{\left(x^2+\pi ^2\right)^2}\\
&=&\frac{k(x)}{\left(x^2+\pi ^2\right)^2}
\end{eqnarray*}
and 
$$p(x)= -x^6-\pi^2 x^4 +\pi^2.x^2 (\pi^2+4) -12\pi^4+\pi^6.$$
Here  $h(0)=h(\pi)=0$ and $h(x_0)=0$.

Elementary computation gives $k(2\pi/3)<0$, while with the use of a computer we get 
$k(17\pi/24)>0$. So $k(x)$ has a root $x_0$ between $2\pi/3\approx 2.0944$ and 
$17\pi/24\approx 2.22529$.

Now, by letting $y=x^2$, we get  
$$p(x)=q(y)= -y^3-\pi^2y^2+\pi^2y(\pi^2+4)+\pi^6-12\pi^4.$$
Since $$q'(y)= -3y^2-2\pi^2y+\pi^2(\pi^2+4),$$ and 
$y>0$, the only root of $q'(y)=0$  is  
$y^*=(2\pi\sqrt{3+\pi^2+3}-\pi^2)/3$, which lies between $\pi$ and $\pi^2$.
One can verify that $q(0)<0$ and $q(\pi)>0$. As $q'(y)>0$ for $y\in(0, y^*)$ 
and $<0$ in $y\in(y^*, \pi^2)$, and as $q(\pi^2)<0$, we get the following:
$q$ is increasing from $q(0)$ to $q(y^*)$ and decreasing from $q(y^*)$ to $q(\pi^2)$. Thus there exist only two roots in $(0, \pi^2 )$ to $q(y)$, let them $z_1$  and $z_2$. 
Clearly, $z_1$ is in $(0, \pi)$ and $y_2$ in $(y^*, \pi^2)$. More precisely, as $q(6)<0$
  and $q(7)>0$, one finds that $y_2 >6$.
These imply that $q(y)<0$ in $(0, z_1)$, $>0$  in $(z_1, z_2)$ and $<$ in $(z_2, \pi^2)$.
In terminology of $P(x)$, we get that
$p(x)<0$ in $(0, \sqrt{z_1})$; $>0$ in $(\sqrt{z_1}, \sqrt{z_2})$, and $<0$ in 
$(\sqrt{z_2}, \pi)$.
In  $(0, \sqrt{z_1})$ clearly $h'(x)<0$, so $h(x)<h(0)=0$; similarly in $(z_2, \pi)$ one 
has $h'(x)<0$, so $h(x)>h(\pi)=0$.
Remains the interval $(\sqrt{z_1}, \sqrt{z_2})$.  As $z_1<\pi$ and  $z_2>6$, we get that  
$$\sqrt{z_1}<\sqrt{\pi}<2< x_0 < \sqrt{6}\approx 2.44949<\sqrt{z_2},$$ 
so we find that $x_0$ lies between $\sqrt{z_1}$ and $\sqrt{z_2}$.
 Then clearly $h(x)<h(x_0)=0$ in $(\sqrt{z_1}, x_0)$, and $h(x)>h(x_0)=0$ in 
$(x_0, \sqrt{z_2})$, so all is done.

Thus the minimum point of $g$ is at $x_0$. As $g(x)$ tends to $1$ when $x$ tends to $0$
 or $\pi$, thus the Redheffer's inequality follows.
On the other hand, we get also $g(x)\geq  g(x_0)$, i.e. the best possible converse to Redheffer's inequality.
Now, with the aid of a computer one can find the more precise approximation 
$x_0\approx 2.12266,$ 
giving  $g(x_0)\approx 0.93012=1/c_1$ so the converse to the Redheffer's inequality holds true.
\end{proof}

\begin{lemma}\label{lema1704} For $a=2.175$, the function
$$Q(x)= -3ax^{2a+2}-2a\pi^ax^{a+2}+2a(a-1)(2a-1)\pi^a x^a +a\pi^{2a}x^a- 
2a(a+1)(a-2)\pi^{2a}$$
has exactly two roots $y_1$ and $y_2$ in $(0,\pi)$.
\end{lemma}

\begin{proof}
We have $Q(1/2)<0,\, Q(\pi/2)>0,$ and $Q(\pi)<0$, so $Q$ has two roots $y_1$ 
in $(1/2, \pi/2)$, resp. $y_2$ in $(\pi/2, \pi)$. To show that $Q$ has no other zeros, we have to consider  
$Q'(x)= x R(x)$, where 
$$R(x)= -3a(2a+2)x^{2a}-2a(a+2)\pi^a x^a -2a^2(a-a)(2a-1)\pi^a x^{a-2}+2a\pi^{2a}.$$
One has further  $$R'(x)= 2a^2x^{a-3} T(x),$$
here $$T(x)=  -3(2a+2)x^{a+2}-(a+2)\pi^a x^2-2a^2(a-1)(a-2)(2a-1)\pi^a.$$
Since $a>2$, we get  $T(x)<0$, so $R'(x)<0.$
One has $Q'(x)= x R(x)$, where $R'(x)<0.$ Since $R(0)>0$ and $R(\pi)<0$, and 
$R(x)$ is strictly decreasing, $R(x)=0$ can have exactly one root $r$ in $(0, \pi)$. Therefore, $Q(x)$ has exactly one extremal point. Since $Q(0)<0,\, Q(\pi)<0$ and $Q$ takes also positive values, clearly $Q(r)$ will be a maximum of $Q(x)$. This shows that $Q$ has exactly two roots in $(0, \pi)$: one in $(0,r)$ and the other one in $(r, \pi)$.

\end{proof}

\begin{theorem} For $x\in(0,\pi)$, the following inequality holds
$$\frac{\sin(x)}{x}<\frac{\pi^a-x^a}{\pi^a+x^a},$$
where $a=87/40=2.175$.
\end{theorem}

\begin{proof}
Inequality can be written as $g(x)=f(x)/\sin (x)>1,$ where 
$$f(x)= \frac{x(\pi^a-x^a)}{\pi^a+x^a}. $$
First of all, similarly to the proof of Theorem \ref{thm16}, one has
$$f’ (x)= \frac{\pi^{2a}-2a\pi^ax^a-x^{2a}}{(\pi^a+x^a)^2},$$   
so  we get  $$h(x)=\frac{K(x)}{(\pi^a+x^a)^2},$$  
where $$K(x)=  (\pi^{2}-2a\pi^a x^a-x^{2a})\sin(x)-x(\pi^{2a}-x^{2a})\cos(x).$$
One finds  $$F(x)=\frac{xP(x)}{(\pi^a+x^a)^3},$$
Where  
$$P(x)=  -x^{3a}-\pi^a x^{2a} +2a(a-1)\pi^a x^{2a-2}+ 
\pi^{2a} x^a-2a(a+1)\pi^{2a}x^{a-2}+\pi^{3a}.$$

Now the proof of Theorem runs as follows: Since $K(1/2)>0$ and $K(\pi/4)<0,\, K(\pi/2)>0$
 and $K(2\pi/3)<0$, we get $x_1$ in $(1/2, \pi/4)$ such that $k(x_1)=0$ and 
$x_2$ in $(\pi/2, 2\pi/3)$ such that $K(x_2)=0$. 
As $P(0)>0$ and $P(\pi/4)<0$. One has 
$$P'(x)=x Q(x),$$
where $Q(x)$ is as in Lemma \ref{lema1704}.
It follows from Lemma \ref{lema1704} that
that, $Q(x)<0$  for $x$ in $(0, y_1)$ and $(y_2, \pi)$, and $Q(x)>0$ for $x$ in $(y_1,y_2)$. This shows that $P(x)$ is strictly decreasing in $(0, y_1)$ and $(y_2, \pi)$ and strictly increasing in $(p_1, p_2)$. This implies that $P(x)$ has a unique root  $p_1$ in $(0, y_1)$, as well as a unique $p_2$ in $(y_1, y_2)$ and $p_3$ in $(y_2, \pi)$. By approximate computation we can see that  $p_1<2/3,\,  p_2>1$ and $p_3>2$,  and $p_1<p_2<p_3$. This shows that  for the roots $x_1$ and $x_2$ of function $k(x)$ one has that $x_1$ is in $(p_1, p_2)$ and $x_2$ in $(p_2,p_3)$.
As $h(0)=h(x_1)=h(x_2)=h(\pi)=0$, we get the following;
\begin{enumerate}
\item for $x \in(0, p_1)$  one has  $h(x)>h(0)=0,$
\item for $x \in(p_1, x_1)$ one has $h(x)>h(x_1)=0,$
\item for $x \in (x_1, p_2)$ one has  $h(x) <h(x_1)=0,$
\item	for $x \in (p_2, x_2)$ one has  $h(x)< h(x_2)=0,$
\item	for $x \in (x_2, p_3)$ one has  $h(x)>h(x_2)=0,$
\item	for $x \in (p_3, \pi)$ one has  $h(x)>h(\pi)=0.$
\end{enumerate}
From the above it follows that  $x_1$ is a local minimum point , while $x_2$ a local maximum point of $g(x)$. 
Clearly, $g(x_1)>\lim g(x)$, when $x$ tends to zero, $=1$  and $\lim g(x)$ at $x =\pi$ is  
$a/2>1.$ This completes the proof.

\end{proof}

\section{A hyperbolic analogue }

\begin{lemma}\label{lem2} For $x\in(0,\infty)$,
$$x+ \sinh(x)>x^2\coth(x/2).$$
\end{lemma} Let $h_2(x)=x+ \sinh(x)-x^2\coth(x/2)$.
\begin{proof} One has
$$h_2'(x)=\frac{(\sinh(x)-x)^2}{2\sinh(x/2)^2}>0,$$
and $h_2(x)>\lim_{x\to 0}h_2(x)=0$. Thus, inequality holds.
\end{proof}

%

\begin{theorem}\label{cor2} For $x\in(0,\infty)$, we have
\begin{equation}\label{redhyp}
\frac{\sinh(x)+x}{\sinh(x)-x}>\frac{12}{x^2}.
\end{equation}
\end{theorem}

\begin{proof} Let 
$$f_2(x)=\frac{x^2(\sinh(x)+x)}{\sinh(x)-x},$$
for $x\in(0,\infty)$. We get
$$\frac{(\sinh(x)-x)^2}{2x}f_2'(x)=(\sinh(x/2))((x+\sinh(x))-x^2\cosh(x/2)^2),$$
which is positive by Lemma \ref{lem2}.
Thus, $f_2$ is strictly increasing in $x\in(0,\infty)$, and the inequality \eqref{redhyp}
follows from $f_2(x)>\lim_{x\to 0}f_2(x)=12$.
Therefore, we have an analogue of the second inequality in the circular case, and this is \eqref{redhyp}.
When $x^2<12$, then \eqref{redhyp} becomes
$$\frac{\sinh(x)}{x}<\frac{12+x^2}{12-x^2}.$$
\end{proof}

It is interesting to observe that
$$\frac{12+x^2}{12-x^2}<\frac{\pi^2+x^2}{\pi^2-x^2},$$
if $x<\pi$. Indeed this becomes: $\pi^2<12$.




\font\fFt=eusm10 
\font\fFa=eusm7  
\font\fFp=eusm5  
\def\K{\mathchoice
{\hbox{\,\fFt K}}
{\hbox{\,\fFt K}}
{\hbox{\,\fFa K}}
{\hbox{\,\fFp K}}}
\def\T{\mathchoice
{\hbox{\,\fFt T}}
{\hbox{\,\fFt T}}
{\hbox{\,\fFa T}}
{\hbox{\,\fFp T}}}


\begin{abstract}
Here authors establish the sharp inequalities for classical beta function by studying the inequalities of trigonometric sine function. 
\end{abstract}

\maketitle

%

\section{Introduction}

For $x,y>0$,  the classical \emph{gamma function} $\Gamma$, the {\it digamma function} $\psi$ and the \emph{beta function}  $B(\cdot,\cdot)$ are defined by
   $$\Gamma(x) = \int^\infty_0 e^{-t}t^{x-1}\,dt,\quad
       \psi(x) = \frac{\Gamma'(x)}{\Gamma(x)},\quad
        B(x,y) = \frac{\Gamma(x)\Gamma(y)}{\Gamma(x+y)},$$
respectively. The study of these functions has become crucial because of its dynamic applications in the field of various branches of engineering and mathematics \cite{roy}.
Since last half century numerous authors has given the several functional inequalities of these function by using different approaches, e.g., see \cite{alz2,alz3,dab,qv}. In this paper, we establish the inequalities for beta function by studying the well-known Jordan's inequality \cite{kvv,mit,neusan}.

The functions $\Gamma$ and $\psi$ satisfy the following recurrence relation
\begin{equation}\label{recgam}
\Gamma(1+x)=x\Gamma(x),\quad \psi(1+x)=x\psi(x).
\end{equation}

 Weierstrass gave the following infinite production definitions of gamma function and
the sine function
$$\frac{1}{\Gamma(x)}=xe^{\gamma x}\prod_{n=1}^\infty\left(1+\frac{x}{n}\right)e^{-x/e},\quad
\sin(\pi x)=\pi x\prod_{n\neq 0}^\infty\left(1-\frac{x}{n}\right)e^{x/n},$$
where $\gamma$ is the Euler-Mascheroni constant \cite{as} defined by
$$\gamma=\lim_{x\to \infty}\left(\sum_{k=1}^\infty \frac{1}{k}-\log(n)\right)\approx 0.57721.$$
These definitions gives the following relation
\begin{equation}\label{eulref}
\Gamma(t)\Gamma(1-t)=\frac{\pi}{\sin (\pi t)},\quad t\notin \mathbb{Z},
\end{equation}
which is known as the Euler's reflection formula \cite[6.1.17]{as}.
We refer to reader to see \cite{dab} for the historical background and the properties of gamma and beta function.

Dragomir et al. \cite{dab} established the following inequality
\begin{equation}\label{ineq0704}
B(x,y)\leq \frac{1}{xy},\quad x,y\in(0,1),
\end{equation} 
which was refined by Alzer \cite{alz1} as follows
\begin{equation}\label{ineq0704a}
\frac{1}{xy}\left(1-a\frac{1-x}{1+x}\frac{1-y}{1+y}\right)< B(x,y)<
\frac{1}{xy}\left(1-b\frac{1-x}{1+x}\frac{1-y}{1+y}\right), ,\quad x,y\in(0,1),
\end{equation}
with the best possible constants $a=2\pi^2/3-4\approx 2.57973$ and $b=1$.
Recently, the second inequality in \eqref{ineq0704} was refined by Iv\'ady\cite {ivady}
\begin{equation}\label{ineq0704b}
\frac{1}{xy}\left(x+y-xy\right)\leq B(x,y)\leq
\frac{1}{xy}\frac{x+y}{1+xy},\quad x,y\in(0,1).
\end{equation}

\begin{lemma}\label{newlema}
Let $g(x)=f(x)/\sin (x)$ for $x \in (0, \pi)$. Then $\sin(x)^2g'(x)= h(x)$, and  the
sign of $h'(x)$ depends on the sign of $F(x)= f(x)+ f''(x)$.
\end{lemma}

\begin{proof} One has  
$$\sin (x)^2 g'(x)= f'(x)\sin(x)-f(x)\cos (x)=h(x),$$ and 
$$h'(x)= (f(x)+f''(x))\sin(x)=F(x)\sin(x).$$
As $\sin (x)>0$ for all $x \in (0,\pi)$, the result follows.
\end{proof}

\begin{theorem} For $x\in(0,1)$, we have
\begin{equation}\label{ineq1104}
\frac{3 (1-t)}{\pi  t^2-\pi  t+\pi }<\frac{\sin (\pi  t)}{\pi  t}<\frac{\pi 
   (1-t)}{\pi  t^2-\pi  t+\pi },
\end{equation}


\begin{equation}\label{ineq0505}
1-(2-t)t^2<\frac{\sin(\pi t)}{\pi t}<\frac{16}{5\pi}\left(1-(2-t)t^2\right).
\end{equation}

\end{theorem}

\begin{proof} Let  $g(x)= f(x)/\sin (x)$  for $x\in(0, \pi)$,
where $$f(x)= (\pi x-x^2)/(\pi^2-\pi x+x^2).$$
We get 
$$\frac{(\pi^2-\pi x+x^2)^3}{x(\pi-x)}F(x)=(\pi^2-\pi x+x^2)^2- 6\pi^2 = A(x)B(x),$$ 
where
$B(x)>0$ always, and 
$A(x)=x^2-\pi x+\pi^2- \pi\sqrt{6}$.
The roots of equation $A(x)=0$ are  
$x_1=  (\pi-\sqrt{4\cdot 6^{1/2}\pi-\pi^2})/2$  which is in $(0, \pi/2)$, and  
$x_2=  (\pi+ \sqrt{4\cdot 6^{1/2}\pi-\pi^2})/2$  which is in $(\pi/2, \pi)$.
Let $x\in(0, x_1)$, then  $A(x)>0$, so $F(x)>0$, giving $h'(x)>0$. 
This implies $h(x)>h(0)=0$, so $g'(x)>0$ by Lemma \ref{newlema}.
Again, let $x\in[x_1, \pi/2)$, then $A(x)\geq 0$, giving  $F(x)\leq 0$, i.e. $h'(x)\leq 0$. 
This implies $h(x)>h(\pi/2)=0$. So $g'(x)\geq 0$ here too.  We have proved that $g'(x)>0$ for all $x$ in $(0, \pi/2)$.
Let now $x$ in $(\pi/2, x_2)$. Then $A(x)<0$, so  $h'(x)<0$, implying  $h(x)<h(\pi/2)=0$. 
For $x$ in $[x_2, \pi)$ one has  $h'(x)\geq 0$, so $h(x)\leq h(\pi)=0$. Therefore, for all $x$ in $(\pi/2, \pi)$ one has $h(x)<0$, i.e. $g'(x)<0$ here. In both cases we had  $g'(x)=0$ only for $x=\pi/2$. Consequently, the function $g$ is strictly increasing in $(0,\pi/2)$ and 
strictly decreasing decreasing $(\pi/2,\pi)$, and attains maximum $1/3$ at $x=\pi/2$ as well
as $g$ tends to $1/\pi$ when $x$ tends to $0$ or $\pi$. This implies the proof of 
\eqref{ineq1104} if we let $x=\pi t$.

For the proof of \eqref{ineq0505}, let 
$$f(x)=\frac{\sin(x)}{x(x^3-2\pi x^2+\pi^3)}.$$

A simple calculation gives

$$(x \,k(x))^2 f'(x)= g(x),$$ 
where $k(x)= x^3-2\pi x^2 +\pi^3$ and
$$g(x)=\cos (x).(x^4-2\pi x^3 +\pi^3x)- \sin (x).(4x^3 -6\pi x^2 +\pi^3).$$
It is immediate that  $g(0)= g(\pi/2)= g(\pi)=0$.
One has $g'(x)= -x\sin(x)h(x),$
with
$$h(x)= x^3-2\pi x^2+12x+\pi^3-12\pi.$$

Here $h(0)=  \pi(\pi^2-12)<0$, $h(\pi/2)= 5\pi^3/8-6\pi >0$  as  $5\pi^2 >48$ and $h(\pi)=0.$
Further  $h'(x)= 3x^2-4\pi x +12,$  and $h''(x)= 2(3x-2\pi)$. Here  $\pi/2< 2\pi/3< \pi$.
The roots of $h'(x)=0$  are $x_1= (2\pi-2\sqrt{\pi^2-9})/3\approx 1.47271$, which is in $(0, \pi/2)$, 
and  $x_2=(2\pi+2\sqrt{\pi^2-9})/3\approx 2.71608$, which is in $(2\pi/3, \pi).$
Therefore, $h(x)$ is strictly increasing in $(0, x_1)$, and $(x_2, \pi)$, while strictly decreasing in $(x_1, x_2). $

Let $x \in (0, \pi/2)$, then as $h(0)<0$, $h(\pi/2)>0$, $h$ has a single root $x_0$, and a maximum point in $x_1$. Thus $h(x)<0$ in $(0, x_0)$, and $h(x)>0$ in $(x_0, \pi/2)$. Therefore, $g'(x)>0$ for $x\in(0, x_0)$ and $g'(x)<0$ in $(x_0, \pi/2)$. Thus 
$g(x)>g(0)=0$ in $(0, x_0)$ and 
$g(x)>g(\pi/2)=0$  $(x_0, \pi/2)$. In all cases,  $g(x)>0$ for $x \in (0, \pi/2)$. This means that, $f(x)$ is strictly increasing in $(0, \pi/2)$.

When $x$ is in $(\pi/2, \pi)$, the proof runs as above, by remarking that by $h(2\pi/3)<0$, there exists a unique $x^*_0 \in (\pi/2, \pi)$ such that $h(x^*_0)=0$. Since $h(x)>0$  in 
$(x^*_0, \pi)$ and $h(x)<0$ in $(x^*_0, \pi)$  we get that $g(x)<g(\pi/2)=0$ in 
$(\pi/2, x^*_0)$, while  $g(x)<g(\pi)=0$ in  $(x^*_0, \pi)$, so in all case $g(x)<0$, 
when $x$ is in $(\pi/2, \pi)$. Thus $f(x)$ is strictly decreasing in $(\pi/2, \pi)$.
This completes the proof.
\end{proof}

%

The inequalities in \eqref{ineq1104} and \eqref{ineq0505} are not comparable.
From the proof of \eqref{ineq0505}
we get the following corollary.
\begin{corollary} For $x\in(0,\pi/2)$, we have
$$\frac{x^3-2 \pi  x^2+\pi ^3}{4 x^3-6 \pi  x^2+\pi ^3}>\frac{\tan(x)}{x},$$
inequality reverses for $x\in(\pi/2,\pi)$.
\end{corollary}

%
%


\begin{theorem}\label{mainthm} For we have
\begin{enumerate}
\item $\displaystyle\frac{\alpha}{xy}\frac{x+y}{1+xy}<B(x,y)<
\displaystyle\frac{\beta}{xy}\frac{x+y}{1+xy},\quad x\in(0,1)\,
{\rm with}\, y=1-x,$\\
with the best possible constants $\alpha=5\pi/16\approx 0.98175$ and $\beta=1$,
\item $B(x,y)<\displaystyle\frac{1}{xy}\frac{x+y}{1+xy},\quad x,y\in(0,1)$,\\
inequality reverses for $x>1$.
\end{enumerate}
\end{theorem}

\begin{proof} Utilizing \eqref{eulref}, the first inequality in \eqref{ineq0505} can be written as 
$$t(1-t)(1+t(1-t))<\frac{1}{\Gamma(t)\Gamma(1-t)},$$
which is equivalent to 
$$\frac{t(1-t)(1+t(1-t))}{t+1-t}<\frac{\Gamma(t+1-t)}{\Gamma(t)\Gamma(1-t)}.$$
Letting $x=t$ and $y=1-t$, we get the first inequality. The second inequality follows similarly from the second inequality of \eqref{ineq0505}. This completes the proof.
\end{proof}

\begin{remark}
The inequality  
$$1-\frac{z}{\pi}<\frac{\sin(z)}{z},\quad z\in(0,\pi),$$
can be written as 
$$\Gamma\left(1+\frac{z}{\pi}\right)\Gamma\left(1-\frac{z}{\pi}\right)<\frac{1}{1-z/\pi},$$
by \eqref{eulref}.
This implies \eqref{ineq0704} if we let $x=z/\pi$ and $y=1-z/\pi$.
\end{remark}

%
%
%
\begin{lemma}\label{lema1004} We have
\begin{enumerate}
\item $\psi(1+x)-\psi(x+y)<\displaystyle\frac{1-y}{x+y-xy},\quad x>1, \,y\in(0,1)$,\\
\item $\psi(2-x)-\psi(1+x)<\displaystyle\frac{1-2x}{1-(1-x)x},\quad x\in(0,1/2)$,\\
inequality reverses for $x\in(1/2,1)$.
\end{enumerate}
\end{lemma}

\begin{proof} For $x>1$ and $y\in(0,1)$, we define
$$g_x(y)=\psi(1+x)-\psi(x+y)-\displaystyle\frac{1-y}{x+y-xy}.$$
Differentiating with respect to $y$ we get
\begin{eqnarray*}
g_x''(y)&=&-\frac{2 (1-x)^2 (1-y)}{(x (1-y)+y)^3}-\frac{2 (1-x)}{(x (-y)+x+y)^2}-
\psi''(x+y)\\
&=&\frac{2 x-2}{(x (1-y)+y)^3}-\psi''(x+y)>0,
\end{eqnarray*}
since $\psi''(x+y)<0$. Thus, $g_x$ is convex in $y$, clearly $g_x(0)=f_x(1)=0$. This implies that the function $g_x$ lies under the line segment joining origin and the point $(1,0)$. Hence the proof is obvious now.

For (2),  write
$$f(x)=\psi(2-x)-\psi(1+x)-\frac{1-2x}{1-(1-x)x},\quad x\in(0,1).$$
One has 
\begin{eqnarray*}
f''(x)&=&\left(\frac{2 (2 x-1)^2}{(1-(1-x) x)^3}-\frac{2}{(1-(1-x) x)^2}\right) (2 x-1)\\
& &-\frac{4 (2 x-1)}{(1-(1-x) x)^2}+\psi''(2-x)-\psi''(x+1)\\
&=&\frac{2 (x-2) (x+1) (2 x-1)}{((x-1) x+1)^3}+\psi''(2-x)-\psi''(x+1).
\end{eqnarray*}
Clearly, the function $\psi''$ is increasing and negative. So, it is not difficult to see that
$f''$ is positive for $x\in(0,1/2)$, and negative for $x\in(1/2,1)$. This implies the convexity and concavity of $f$ in $x\in(0,1/2)$ and $x\in(1/2,1)$, respectively. Clearly, 
$f(0)=f(1/2)=f(1)=0$. This completes the proof. 
\end{proof}

\begin{theorem}\label{thm1004} We have
\begin{enumerate}
\item $\displaystyle \frac{1}{xy}(x+y-xy)>B(x,y),
\quad x>1,\,y\in(0,1)$,\\
inequality reverses for $x\in(0,1)$,
\item $B(x,y)<\displaystyle\displaystyle\frac{\pi}{3}\frac{1}{xy}(x+y-xy),
\quad x\in(0,1),\,{\rm with}\, y=1-x$.
\end{enumerate}
\end{theorem}

\begin{proof} The inequality in (1) can be written as
$$h_y(x)=\log(\Gamma(1+x))+\log(\Gamma(1+y))-\log(\Gamma(x+y))+\log(x+y-xy)>0.$$
Clearly, $h_y(1)=0$.
Differentiation with respect to $x$ yields
$$h_y'(x)=\psi(1+x)-\psi(x+y)-\displaystyle\frac{1-y}{x+y-xy}=g_x(y),$$
which is negative by Lemma \ref{lema1004}(1). thus the function $h_y(x)$ is decreasing in $x>1$, this implies (1).
For part (2),
let $$h(x)= \log \left(\frac{1}{3} \pi  (1-(1-x) x)\right)-
\log (\Gamma (2-x))-\log (\Gamma (x+1)),$$
clearly $h(1/2)=0$.
One has,
$$h'(x)=\psi(2-x)-\psi(1+x)-\frac{1-2x}{1-(1-x)x},$$
which is positive in $x\in(0,1/2)$ and negative in $x\in(1/2,1)$ by Lemma \ref{lema1004} (2). This implies that
$h$ is decreasing in $x\in(0,1/2)$ and increasing in $x\in(1/2,1)$. Thus the proof follows.
\end{proof}
